\documentclass[a4paper,leqno,12pt]{amsart}
\usepackage[utf8]{inputenc}

\usepackage[top=25mm, bottom=25mm, left=20mm, right=20mm]{geometry}
\usepackage{amsmath, amstext, amssymb, enumitem, mathtools}
\usepackage{caption}

\usepackage{hyperref}
\usepackage{nameref}

\hypersetup{colorlinks = true, urlcolor = blue, linkcolor = blue, citecolor = red}

\newcommand{\NN}{\mathbb{N}}
\newcommand{\RR}{\mathbb{R}}
\newcommand{\CC}{\mathbb{C}}
\newcommand{\BB}{\mathbb{B}}
\newcommand{\VV}{\mathbb{V}}

\newcommand{\me}{\mathrm{e}}
\newcommand{\dint}{\text{\normalfont{d}}}

\DeclarePairedDelimiter{\norm}{\lVert}{\rVert}
\DeclarePairedDelimiter{\abs}{\lvert}{\rvert}
\DeclarePairedDelimiter{\inner}{\langle}{\rangle}

\numberwithin{equation}{section}

\theoremstyle{plain}
\newtheorem{theorem}[equation]{Theorem}
\newtheorem{lemma}[equation]{Lemma}
\newtheorem{proposition}[equation]{Proposition}

\theoremstyle{definition}
\newtheorem{remark}[equation]{Remark}

\begin{document}

\title[Jacobi heat kernel estimates on a ball and on a simplex for all parameters]{Sharp estimates for Jacobi heat kernels \\ on a ball and on a simplex \\ in the full range of parameters}

\author{Dawid Hanrahan}
\address{Dawid Hanrahan (\textnormal{dawid.hanrahan@pwr.edu.pl})
\newline WUST -- Wroc\l aw University of Science and Technology, 50-370 Wroc\l aw, Poland}

\begin{abstract}
We establish genuinely sharp heat-kernel estimates for orthogonal
expansions on the Euclidean unit ball and on the unit simplex throughout
their natural parameter ranges. For the ball weight
$W_\mu(x)=(1-\norm{x}^2)^{\mu-1/2}$, we extend the known estimate from the range
$\mu\geq0$ to every $\mu > -1/2$. For the simplex, we prove the corresponding
estimate for every $\kappa\in(-1/2,\infty)^{d+1}$. The estimates retain
exactly the form known in the classical ranges. The proof is based on
analytic continuation of the classical product and reduction formulas.
Below the integrability threshold, the ball formula decomposes into an
endpoint term and a derivative correction, while the simplex formula
produces endpoint and mixed-derivative terms indexed by the negative
coordinates. In the lowest simplex parameter range we use the auxiliary
even ultraspherical kernel introduced by Nowak, Sj\"ogren, and Szarek.


\end{abstract}

\maketitle

\section{Introduction}

Genuinely sharp heat-kernel bounds describe not only the Gaussian decay
but also the precise polynomial factors arising from the geometry and
the underlying measure. Nowak, Sj\"ogren, and Szarek obtained such
bounds for Jacobi expansions, compact rank-one symmetric spaces, the
unit ball, and the unit simplex in \cite{NSS19}. In the ball and simplex
settings, their proofs rely on product formulas that reduce the
multidimensional kernels to one-point Jacobi heat kernels. The reduction strategy developed in \cite{NSS19} has also been applied
to obtain genuinely sharp heat-kernel estimates on conic and
double-conic domains; see \cite{HK23,Han24}. These
representations are particularly effective, but in their classical form
they require all auxiliary measures to be finite.

Further developments include sharp heat-kernel estimates for
Fourier--Dini expansions \cite{LN24} and a unified treatment of sharp
Jacobi heat-kernel bounds with quantitative control of the
multiplicative constants \cite{NPS25}.

For the weighted unit ball, the natural parameter range is
$\mu>-1/2$, whereas the sharp estimate in \cite[Theorem~5.1]{NSS19}
was proved for $\mu\geq0$. The restriction comes from the measure
$\dint\Pi_{\mu-1/2}$ in the reduction formula. When
$-1/2<\mu<0$, its formal density is not integrable at the endpoints,
although its primitive remains integrable. Thus, the obstruction lies in
the representation rather than in the ball operator or its heat
semigroup.

The same phenomenon occurs on the simplex, where the weight has parameters
$\kappa_j>-1/2$ for $1\leq j\leq d+1$. The classical simplex formula contains
the tensor product of the measures
$\dint\Pi_{\kappa_j-1/2}$ and was used in \cite[Theorem~6.1]{NSS19}
under the assumption $\kappa_j\geq0$. If one or more coordinates of
$\kappa$ are negative, the corresponding integrations must be
regularized simultaneously.

Our argument follows the continuation method developed by Nowak, Sj\"ogren, and Szarek in \cite{NSS24}.
That method extends the Dijksma--Koornwinder formula and the Jacobi
heat-kernel reduction below the classical threshold. At the polynomial
level, endpoint values are subtracted from the integrands before the
parameters are continued. Integration by parts then replaces each
singular pairing by an endpoint contribution and an ordinary integral
involving a derivative. After summation in the degree, this gives a
regularized heat-kernel reduction formula.

On the ball there is one singular auxiliary variable. Consequently, the
continued kernel is the sum of one endpoint term and one derivative
correction, and these terms provide the lower estimate in complementary
geometric regimes. The simplex has a parallel but genuinely
multi-parameter structure. Denoting
$J_< \coloneqq \{j:\kappa_j<0\}$, one obtains a reduction formula with $2^{\abs{J_<}}$
non-negative terms indexed by subsets of $J_<$. A coordinate contributes
an endpoint value or a derivative according to whether
$\sqrt{x_jy_j}$ is smaller or larger than the heat scale $t$. A further
difference is that the inner ultraspherical parameter may lie in
$(-3/2,-1]$. In that range an ordinary Jacobi heat kernel is no longer
available, and we use the auxiliary even ultraspherical kernel from
\cite[Section~4]{NSS24}.

The resulting estimates have the same form as their classical
counterparts. In particular, the conditions $\mu\geq0$ on the ball and $\kappa_j\geq0$ for each coordinate on the simplex come from the positive integral representations, not from changes in the heat-kernel geometry.

Our main results are stated below. The notation used in their statements is introduced in Section~\ref{sec:preliminaries}, and the
holomorphy condition required for the continuation argument is discussed in 
Appendix~\ref{app:holomorphy}.

\begin{theorem}[Sharp estimate on the ball in the full range]
\label{thm:main-ball}
Fix $d\geq 2$, $\mu > -1/2$, and $T>0$. Then one has
\begin{equation}
\label{eq:main-estimate}
    h_t^\mu(x,y) \simeq
    [t+\pi-d_{\BB}(x,y)]^{-\lambda_\mu}
    \bigg(t+\frac{\delta(x)\delta(y)}{\pi-d_{\BB}(x,y)}\bigg)^{-\mu}
    t^{-d/2}
    \exp\bigg(-\frac{d_{\BB}(x,y)^2}{4t}\bigg),
\end{equation}
uniformly in $x,y\in \BB^d$ and $0<t\leq T$. Moreover, one has 
\[
    h_t^\mu(x,y)\simeq 1,
\]
uniformly in $x,y\in\BB^d$ and $t\geq T$.
\end{theorem}

\begin{theorem}[Sharp estimate on the simplex in the full range]
\label{thm:main-simplex}
Fix $d\geq2$, $\kappa\in(-1/2,\infty)^{d+1}$, and $T>0$. Then, denoting $x_{d+1} \coloneqq 1-x_1 - \dots - x_d$ and $y_{d+1} \coloneqq 1- y_1 - \dots - y_d$, one has
\begin{equation}
\label{eq:main-simplex-estimate}
    H_t^\kappa(x,y)
    \simeq
    \bigg(\prod_{j=1}^{d+1}
    \big(t+\sqrt{x_jy_j}\big)^{-\kappa_j}\bigg)
    t^{-d/2}
    \exp\bigg(-\frac{d_{\VV}(x,y)^2}{t}\bigg),
\end{equation}
uniformly in $x,y\in\VV^d$ and $0<t\leq T$. Moreover, one has
\[
    H_t^\kappa(x,y)\simeq1,
\]
uniformly in $x,y\in\VV^d$ and $t \geq T$.
\end{theorem}

As usual, in \eqref{eq:main-estimate} the quotient
\[
    \frac{\delta(x)\delta(y)}{\pi-d_{\BB}(x,y)}
\]
is understood by continuous extension at antipodal boundary points; see \eqref{eq:delta-A-B} and \eqref{eq:ball-distance}.

The paper is organized as follows. Section~\ref{sec:preliminaries}
collects the common Jacobi, ball, and simplex preliminaries.
Section~\ref{sec:ball-case} treats the ball, from the regularized
reproducing formula through the proof of the sharp bounds.
Section~\ref{sec:simplex-case} develops the corresponding
multi-parameter argument on the simplex. Finally, Appendix~\ref{app:holomorphy}
discusses the required holomorphy condition. 

\section{Preliminaries}
\label{sec:preliminaries}

Throughout the paper, we use standard notation. We write $X\lesssim Y$ if $X\leq CY$ with a positive constant $C$ independent of all significant quantities. We write $X\simeq Y$ if simultaneously $X\lesssim Y$ and $Y\lesssim X$. The implicit constants $C$ may change from line to line. Unless explicitly stated otherwise, they may depend on fixed structural parameters, such as $d$ and $T$, and the parameters of the measures or kernels under consideration, but not on the running space or time variables.

The minimum and maximum of two quantities will be denoted by $\wedge$ and $\vee$, respectively.

\subsection{Measures}
For $\alpha>-1/2$, let $\dint \Pi_\alpha$ be the probability measure on $[-1,1]$ given by
\begin{equation}
\label{eq:Pi-positive}
    \dint \Pi_\alpha(u)
    \coloneqq \frac{\Gamma(\alpha+1)}{\sqrt{\pi}\Gamma(\alpha+1/2)}(1-u^2)^{\alpha-1/2} \, \dint u.
\end{equation}
For $\alpha=-1/2$, we consider the weak limit of the measures in \eqref{eq:Pi-positive}, that is,
\[
    \dint \Pi_{-1/2} \coloneqq \frac{\delta_{-1}+\delta_1}{2},
\]
where $\delta_x$ denotes the Dirac measure at $x$.
For $-1<\alpha<-1/2$, we use the analytic-continuation convention from the Jacobi setting \cite{NSS24}. Specifically, for each $\alpha \in \CC$ such that $\operatorname{Re} \alpha > -1$ and $\alpha \neq -1/2$, let $\Pi_\alpha$ be the continuous odd function on $(-1,1)$ defined on $(0,1)$ by
\begin{equation*}
    \Pi_\alpha(u)
    \coloneqq \frac{\Gamma(\alpha+1)}{\sqrt{\pi}\Gamma(\alpha+1/2)}
    \int_0^u (1-w^2)^{\alpha-1/2} \, \dint w.
\end{equation*}
For real $\alpha > -1/2$, its derivative coincides with the density of $\dint \Pi_\alpha$, which makes this notation consistent. On the other hand, for arbitrary $\alpha \in \CC$ such that $\operatorname{Re} \alpha > -1$ and $\alpha \neq -1/2$, we can consider the distribution derivative $\dint \Pi_\alpha$, which is a local complex measure in $(-1,1)$. In particular, for real $-1<\alpha<-1/2$, this distribution derivative has a density that is negative and not integrable in $(-1,1)$, but $\Pi_\alpha$ itself is a density of a finite measure on $[-1,1]$ because
\begin{equation}
\label{eq:Pi-size}
    \abs{\Pi_\alpha(u)} \simeq_\alpha |u| (1-|u|)^{\alpha+1/2}, \qquad u \in (-1,1),
\end{equation}
as stated in \cite[Lemma~2.2]{NSS15} or \cite[Lemma~2.1]{NSS24}.

We shall use the following one-variable regularization principle, cf.~\cite[Proposition~2.2]{NSS24}.

\begin{proposition}
\label{prop:one-variable-Pi}
Consider the open set
\[
    H_\alpha \coloneqq \{ \alpha \in\CC: \operatorname{Re} \alpha > -1 \}.
\]
Assume that $\varphi_\alpha(u)$ is continuous in $(u,\alpha)\in[0,1]\times H_\alpha$ and holomorphic in $\alpha\in H_\alpha$ for each fixed $u\in(0,1)$. If, for every compact subset $K\subset H_\alpha$, there exists a constant $C_K \in (0,\infty)$ such that
\[
    \abs{\varphi_\alpha(u)}\leq C_K(1-u), \qquad (u,\alpha) \in [0,1] \times K,
\]
then the function
\[
    I(\alpha) \coloneqq \int_{(0,1]}\varphi_\alpha(u) \, \dint\Pi_{\alpha}(u),
    \qquad \alpha \in H_\alpha,
\]
is well defined and holomorphic in $H_\alpha$.
Moreover, $I(-1/2) = 0$. Next, consider the open set
\[
    H_{\alpha,\beta} \coloneqq \{ (\alpha, \beta) \in\CC^2 : \operatorname{Re} \alpha > -1 \text{ and } \operatorname{Re} \beta > -1 \}.
\]
Assume that $\varphi_{\alpha,\beta}(u)$ is continuous in $(u,\alpha, \beta)\in[0,1]\times H_{\alpha, \beta}$ and holomorphic in $(\alpha, \beta) \in H_{\alpha,\beta}$ for each fixed $u\in(0,1)$. If, for every compact subset $K\subset H_{\alpha, \beta}$, there exists a constant $C_K \in (0,\infty)$ such that
\[
    \abs{\varphi_{\alpha, \beta}(u)} \leq C_K(1-u), \qquad (u,\alpha,\beta) \in [0,1] \times K,
\]
then the function
\[
    I(\alpha, \beta) \coloneqq \int_{(0,1]}\varphi_{\alpha, \beta}(u) \, \dint\Pi_{\alpha}(u),
    \qquad (\alpha,\beta) \in H_{\alpha, \beta},
\]
is well defined and holomorphic in $H_{\alpha, \beta}$. Moreover, $I(-1/2, \beta) = 0$. Finally, assume that $\varphi_{\alpha,\beta}(u,v)$ is continuous in $(u,v,\alpha, \beta)\in[0,1]^2\times H_{\alpha, \beta}$ and holomorphic in $(\alpha, \beta) \in H_{\alpha,\beta}$ for each fixed $(u,v)\in(0,1)^2$. If, for every compact subset $K\subset H_{\alpha, \beta}$, there exists a constant $C_K \in (0,\infty)$ such that
\[
    \abs{\varphi_{\alpha, \beta}(u,v)} \leq C_K(1-u)(1-v), \qquad (u,v,\alpha,\beta) \in [0,1]^2 \times K,
\]
then the function
\[
    I(\alpha, \beta) \coloneqq \int_{(0,1]^2}\varphi_{\alpha, \beta}(u,v) \, \dint\Pi_{\alpha}(u)
    \, \dint\Pi_{\beta}(v),
    \qquad (\alpha,\beta) \in H_{\alpha, \beta},
\]
is well defined and holomorphic in $H_{\alpha, \beta}$.
Moreover, $I(-1/2, \beta) = I(\alpha, -1/2) = 0$.
\end{proposition}

\begin{proof}
We present a detailed proof of the first assertion; the variants involving the additional parameter $\beta$ and the double integral follow by the same argument. 

Let $K\subset H_\alpha$ be compact. Then $\operatorname{Re} \alpha\geq -1+\varepsilon$ on $K$ for some $\varepsilon>0$, and the function
\[
    \alpha\longmapsto
    \frac{\Gamma(\alpha+1)}{\sqrt{\pi}\Gamma(\alpha+1/2)}
\]
is holomorphic, hence bounded, on $K$. For $\alpha\neq -1/2$ the integral can be written as
\[
    I(\alpha)
    =
    \frac{\Gamma(\alpha+1)}{\sqrt{\pi}\Gamma(\alpha+1/2)}
    \int_0^1
    \varphi_\alpha(u)(1-u^2)^{\alpha-1/2}\,\dint u.
\]
Moreover,
\[
    \abs{\varphi_\alpha(u)(1-u^2)^{\alpha-1/2}}
    \lesssim_K
    (1-u)(1-u^2)^{\operatorname{Re} \alpha-1/2}
    \le
    (1-u)^{\varepsilon-1/2},
\]
and the last function is integrable on $(0,1)$. Thus, the integral is locally uniformly convergent in $H_\alpha$ and defines a holomorphic function of $\alpha$ by the standard theorem on parameter-dependent holomorphic integrals.

At $\alpha=-1/2$ the preceding formula has a removable value equal to zero, because $1/\Gamma(\alpha+1/2)$ vanishes at $\alpha=-1/2$ and the integral remains finite. This agrees with the endpoint convention, since the assumed bound implies $\varphi_{-1/2}(1)=0$, and therefore
\[
    \int_{(0,1]}\varphi_{-1/2}(u)\,\dint\Pi_{-1/2}(u)
    =
    \frac12\varphi_{-1/2}(1)
    =
    0.
\]
This completes the proof of the first assertion.

For the version with the passive parameter $\beta$, the same domination is uniform on compact subsets of $H_{\alpha,\beta}$, so the integral is holomorphic in both variables and vanishes when $\alpha=-1/2$. In the double integral case, away from $-1/2 \in \{\alpha,\beta\}$, one writes
\[
    I(\alpha,\beta)
    =
    \frac{\Gamma(\alpha+1)}{\sqrt{\pi}\Gamma(\alpha+1/2)}
    \frac{\Gamma(\beta+1)}{\sqrt{\pi}\Gamma(\beta+1/2)}
    \int_0^1\int_0^1
    \varphi_{\alpha,\beta}(u,v)
    (1-u^2)^{\alpha-1/2}
    (1-v^2)^{\beta-1/2}\,\dint u\,\dint v.
\]
The assumed bound gives the integrable majorant
\[
    (1-u)^{\varepsilon-1/2}(1-v)^{\varepsilon-1/2}
\]
on each compact subset of $H_{\alpha,\beta}$, with some new, possibly smaller, parameter $\varepsilon>0$. Hence, the double integral is holomorphic. Finally, the same removable-singularity argument gives the desired equality
$I(-1/2,\beta)=I(\alpha,-1/2)=0$.
\end{proof}

\subsection{Jacobi and ultraspherical heat kernels}

Fix $\alpha,\beta>-1$. For each $t > 0$, we define the Jacobi heat kernel
\begin{equation}
\label{eq:jacobi-heat-kernel}
    G_t^{\alpha,\beta}(x,y)
    \coloneqq \sum_{n=0}^{\infty}\me^{-tn(n+\alpha+\beta+1)}
    \frac{P_n^{\alpha,\beta}(x)P_n^{\alpha,\beta}(y)}{h_n^{\alpha,\beta}},
    \qquad (x,y)\in[-1,1]^2,
\end{equation}
using the Jacobi polynomials $P_n^{\alpha,\beta}$ and their norms $h_n^{\alpha,\beta}$ in $L^2 \big([-1,1],(1-x)^\alpha(1+x)^\beta \, \dint x \big)$.
When $\alpha,\beta$ are complex with
$\operatorname{Re}\alpha,\operatorname{Re}\beta>-1$, we interpret the
Jacobi polynomials and reciprocal norm coefficients by their explicit
polynomial and gamma-function formulas; see
\cite[Section~2]{NSS24}. The series in
\eqref{eq:jacobi-heat-kernel} then converges locally uniformly in the
parameters and gives their holomorphic extension. No orthogonality
interpretation is attached to this notation for non-real parameters.

In the ultraspherical case, we abbreviate
\[
    G_t^\lambda(x,y) \coloneqq G_t^{\lambda,\lambda}(x,y).
\]
We shall use the sharp one-point estimate
\begin{equation}
\label{eq:one-point-ultra}
    G_t^\lambda(\cos\theta,1)
    \simeq
    t^{-\lambda-1}[t+\pi-\theta]^{-\lambda-1/2}
    \exp\bigg(-\frac{\theta^2}{4t}\bigg),
\end{equation}
which holds uniformly in $(\theta, t) \in [0, \pi] \times (0,T]$ for each fixed $\lambda\geq -1/2$; see~ \cite[Theorem~3.1]{NSS19}.

We shall also need the following derivative identity
\begin{equation}
\label{eq:derivative-ultra}
    \frac{\dint}{\dint z}G_t^\lambda(z,1)
    =2(\lambda+1)\me^{-t(2\lambda+2)}G_t^{\lambda+1}(z,1),
    \qquad z\in(-1,1),
\end{equation}
which holds for each fixed $t > 0$; see~\cite[Lemma~3.8]{NSS19}. 
We note that $z\mapsto G_t^\lambda(z,1)$ is increasing.

\subsection{The ball setting}

We recall the standard ball setting from \cite[Section~5]{NSS19}. Fix a positive integer $d\geq2$. We denote by
\[
    \BB^d \coloneqq \{x\in\RR^d:\norm{x}\leq 1\}
\]
the closed $d$-dimensional Euclidean unit ball.
For each fixed $\mu>-1/2$, we let
\[
    W_\mu(x) \coloneqq (1-\norm{x}^2)^{\mu-1/2},\qquad x\in\BB^d.
\]
We also define the auxiliary parameters
\begin{equation*}
    \lambda_\mu \coloneqq \mu+\frac{d-1}{2}, \qquad \nu_\mu \coloneqq \lambda_\mu-\frac12 = \mu+\frac{d-2}{2}.
\end{equation*}
The associated heat kernel is given by
\begin{equation}
\label{eq:ball-heat-series}
    h_t^\mu(x,y) \coloneqq \sum_{n=0}^{\infty}\me^{-tn(n+2\lambda_\mu)}P_n(W_\mu,x,y), \qquad (x,y) \in \BB^d \times \BB^d,
\end{equation}
where $P_n(W_\mu,x,y)$ is the reproducing kernel of the space of orthogonal polynomials of degree $n$ with respect to the measure $\dint W_\mu$. If $\mu \geq 0$, then $P_n(W_\mu,x,y)$ is given by the formula
\begin{equation}
\label{eq:classical-reproducing-kernel}
    \frac{1}{W_\mu(\BB^d)}\frac{n+\lambda_\mu}{\lambda_\mu}
    \int_{[-1,1]}C_n^{\lambda_\mu}
    \Big(\inner{x,y}+u\sqrt{1-\norm{x}^2}\sqrt{1-\norm{y}^2}\Big)
    \, \dint\Pi_{\mu-1/2}(u),
\end{equation}
where $C_n^\lambda$ is the Gegenbauer polynomial; see
\cite[Theorem~4.3]{Xu}. For real $\mu>-1/2$, the symbols
$P_n(W_\mu,x,y)$ and $h_t^\mu(x,y)$ always denote the genuine
reproducing and heat kernels associated with the positive measure
$\dint W_\mu$. Appendix~\ref{app:holomorphy} constructs holomorphic
extensions of these functions to complex $\mu$ with
$\operatorname{Re}\mu>-1/2$. For non-real $\mu$, the same symbols denote
only those holomorphic extensions; no positive measure or Hilbert-space
interpretation is intended.

We define the auxiliary quantities
\begin{equation}
\label{eq:delta-A-B}
    \delta(x) \coloneqq \sqrt{1-\norm{x}^2},\qquad A \coloneqq \inner{x,y},\qquad B \coloneqq \delta(x)\delta(y).
\end{equation}
Note that $0\leq B\leq1$ and $-1\leq A\leq1-B$. The natural distance in this setting is
\begin{equation}
\label{eq:ball-distance}
    d_{\BB}(x,y) \coloneqq \arccos(A+B)
    =\arccos(\inner{x,y}+\delta(x)\delta(y)).
\end{equation}
Precisely, this distance coincides with the geodesic distance between the lifted points $(x,\delta(x))$ and $(y,\delta(y))$ on the unit sphere in $\RR^{d+1}$. In particular, we see that $0\leq d_{\BB}(x,y)\leq\pi$.

For $\mu\geq0$, the classical reduction formula recalled in \cite[Section~5]{NSS19} reads as follows
\begin{equation}
\label{eq:classical-ball-reduction}
    h_t^\mu(x,y)
    = c_{\mu,d}\int_{[-1,1]}
    G_t^{\nu_\mu} (A+Bu,1) \, \dint\Pi_{\mu-1/2}(u),
\end{equation}
where $c_{\mu,d}$ is the normalizing constant given by
\begin{equation*}
    c_{\mu,d}\coloneqq \frac{2^{2\lambda_\mu}\Gamma(\lambda_\mu+1/2)^2}{\Gamma(2\lambda_\mu+1)W_\mu(\BB^d)}.
\end{equation*}
Notice that if $d\geq2$ and $\mu>-1/2$, then the corresponding ultraspherical parameter satisfies
\begin{equation}
\label{eq:nu-range}
    \nu_\mu=\mu+\frac{d-2}{2}>-\frac12.
\end{equation}
Thus, the inner ultraspherical parameter remains in the standard range $\nu_\mu>-1/2$. Unlike the general Jacobi continuation argument of \cite{NSS24}, the ball case therefore requires regularization only in the outer auxiliary variable.

\subsection{The simplex setting}

We recall the simplex setting from \cite[Section~6]{NSS19}. Let
\[
    \VV^d
    \coloneqq
    \{x\in\RR^d: x_1, \dots, x_d \geq 0 \text{ and } x_1 + \dots + x_d \leq 1 \}.
\]
For
$\kappa=(\kappa_1,\ldots,\kappa_{d+1})\in(-1/2,\infty)^{d+1}$,
we equip $\VV^d$ with the measure
\[
    \dint U_\kappa(x)
    \coloneqq
    \bigg( \prod_{j=1}^{d+1}x_j^{\kappa_j-1/2} \bigg) \,
    \dint x_1\cdots\dint x_d,
\]
where $x_{d+1} \coloneqq 1- x_1 - \dots - x_d$. 
The associated heat kernel has the spectral expansion
\begin{equation}
\label{eq:simplex-heat-series}
    H_t^\kappa(x,y)
    \coloneqq \sum_{n=0}^\infty
    \me^{-tn(n+\lambda_\kappa)}P_n(U_\kappa,x,y),
\end{equation}
where $\lambda_\kappa \coloneqq \abs{\kappa}+(d-1)/2$ and $\abs{\kappa} \coloneqq \kappa_1 + \dots + \kappa_{d+1}$, and
$P_n(U_\kappa,x,y)$ is the reproducing kernel of the space of orthogonal polynomials of degree $n$ with respect to the measure $\dint U_\kappa$.
As before, the parameter $\kappa$ is initially assumed to be real, while complex values of $\kappa$ are introduced later as an auxiliary tool in the analytic continuation argument. The corresponding local holomorphic extensions are
constructed and distinguished from the real kernels in
Appendix~\ref{app:holomorphy}.

For $x,y\in\VV^d$ and $u\in[-1,1]^{d+1}$, set
\begin{equation}
\label{eq:simplex-a-z}
    a_j \coloneqq \sqrt{x_jy_j},
    \qquad
    z(u;x,y) \coloneqq \sum_{j=1}^{d+1}a_ju_j.
\end{equation}
The intrinsic simplex distance is given by
\begin{equation*}
    d_{\VV}(x,y)
    \coloneqq \arccos z(\boldsymbol{1};x,y)
    =\arccos\bigg(\sum_{j=1}^{d+1}\sqrt{x_jy_j}\bigg).
\end{equation*}
In particular, we see that $0\leq d_{\VV}(x,y)\leq\pi/2$. We shall also use the parameters
\begin{equation*}
    \alpha_j \coloneqq \kappa_j-\frac12,
    \qquad
    \rho_\kappa \coloneqq \lambda_\kappa-\frac12
    =\abs{\kappa}+\frac d2-1.
\end{equation*}
Observe that, in the full real range $\kappa \in(-1/2,\infty)^{d+1}$, we have the restriction
\begin{equation}
\label{eq:rho-kappa-range}
    \rho_\kappa>-\frac32.
\end{equation}
In contrast to \eqref{eq:nu-range} from the ball setting, the range in \eqref{eq:rho-kappa-range} goes beyond the range $(-1,\infty)$ of the inner ultraspherical parameter from \eqref{eq:jacobi-heat-kernel} in the case $\alpha = \beta$. This is the reason for introducing the auxiliary kernel
from \cite[Section~4]{NSS24}.

For a fixed even integer $k \geq 0$, we define
\begin{align*} 
    \mathfrak D_k^\rho \coloneqq
    \begin{cases}
(2k+2\rho+1)\Gamma(k+2\rho+1)[2^{2\rho+1}\Gamma(\rho+2)\Gamma(k+\rho+1)]^{-1} & \text{if } k \geq 2, \\
\Gamma(\rho+3/2)[\sqrt\pi\Gamma(\rho+2)]^{-1} & \text{if } k = 0,
\end{cases}
\end{align*}
so that the map $\rho \mapsto \mathfrak D_k^\rho$ is holomorphic in the half-plane $\operatorname{Re} \rho > -3/2$. 
For real $\rho>-3/2$, writing $P^{\rho}_{2n}$ for $P^{\rho, \rho}_{2n}$, we let
\begin{equation}
\label{eq:auxiliary-H-rho}
    \mathcal H_s^\rho(z)
    \coloneqq
    \sum_{n=0}^\infty
    \me^{-s2n(2n+2\rho+1)}
    \mathfrak D_{2n}^\rho P_{2n}^\rho(z).
\end{equation}
For $\rho>-1$, this is the product of $(\rho+1)^{-1}$ and the even part of the ultraspherical heat kernel 
$G_s^\rho(z,1)$. Formula \eqref{eq:auxiliary-H-rho} extends this quantity to $\rho > -3/2$. One also has
\begin{equation}
\label{eq:auxiliary-H-estimate}
    \mathcal H_s^\rho(\cos\theta)
    \simeq
    s^{-\rho-1}\exp\bigg(-\frac{\theta^2}{4s}\bigg),
\end{equation}
uniformly in $0\leq\theta\leq\frac\pi2$ and $0<s\leq T$; see~\cite[Lemmas~5.1 and~5.3]{NSS24}.

Differentiating \eqref{eq:auxiliary-H-rho} repeatedly and using the calculation in \cite[Lemma~4.2]{NSS24}, we obtain
\begin{equation}
\label{eq:auxiliary-H-derivatives}
    \frac{\dint^r}{\dint z^r}\mathcal H_s^\rho(z)
    =2^r\frac{\Gamma(\rho+r+1)}{\Gamma(\rho+2)}
    \me^{-sr(2\rho+r+1)}
    \big[G_s^{\rho+r}(z,1)\big]_{\varepsilon_r}
\end{equation}
for positive integers $r \geq 1$, where $\big[G_s^{\rho+r}(z,1)\big]_{\varepsilon_r}$ denotes the odd part of $G_s^{\rho+r}(z,1)$ if $r$ is odd, and the even part if $r$ is even. Since $\rho+r>-1/2$, all Jacobi kernels on the
right-hand side are in the standard range.

\subsection{The one-variable integral estimate}

The following estimate is the integral engine behind the ball and simplex estimates; cf.~\cite[Lemma~2.1]{NSS19}.

\begin{lemma}
\label{lem:integral-engine}
Fix $\gamma\geq -1/2$ and $\xi\in\RR$. Denoting
\[
    \Phi_{A,B}(w) \coloneqq \arccos(A+Bw), \qquad w \in [-1,1],
\]
we have, uniformly in $0\leq B\leq1$, $-1\leq A\leq1-B$, and $D>0$, the estimate
\begin{align*}
\begin{split}
&\int_{[0,1]}
    \big(\pi-\Phi_{A,B}(w)+D\big)^{-\xi}
    \exp (- \Phi_{A,B}(w)^2 / D )
    \, \dint\Pi_\gamma(w)
\\
&\quad\simeq
    D^{\gamma+1/2}
    \big(\pi-\Phi_{A,B}(1)+D\big)^{-\xi}
    \big(B/(\pi-\Phi_{A,B}(1))+D\big)^{-\gamma-1/2}
    \exp (- \Phi_{A,B}(1)^2 / D ).
\end{split}
\end{align*}
If $B=0$, then the quotient $B/(\pi-\Phi_{A,B}(1))$ is understood as $0$.
\end{lemma}

\section{The ball case}
\label{sec:ball-case}

\subsection{A regularized formula for the reproducing kernel}
\label{sec:regularized-reproducing-kernel}

We first record the polynomial-level counterpart of the continuation
argument. This is the one-variable analogue of the regularized
Dijksma--Koornwinder formula in the Jacobi setting;
cf.~\cite[Sections~3 and~4]{NSS24}. The complex-valued kernel appearing
in the argument is the independent holomorphic extension constructed
in Appendix~\ref{app:holomorphy}; the formula below identifies that
extension and is not being used to define a reproducing kernel for a
complex weight.

For $n\geq0$ and $x,y\in\BB^d$, define the renormalized reproducing kernel
\begin{equation}
\label{eq:renormalized-reproducing-kernel}
    \mathcal P_n^\mu(x,y)
    \coloneqq W_\mu(\BB^d)\frac{\lambda_\mu}{n+\lambda_\mu}P_n(W_\mu,x,y).
\end{equation}
Furthermore, set
\begin{equation*}
    \Theta_n^\mu(u;x,y) \coloneqq C_n^{\lambda_\mu}(A+Bu),\qquad -1\leq u\leq1,
\end{equation*}
and
\begin{equation*}
    \Theta_{n,E}^\mu(u;x,y)
    \coloneqq \frac12\left(\Theta_n^\mu(u;x,y)+\Theta_n^\mu(-u;x,y)\right),
    \qquad 0\leq u\leq1.
\end{equation*}

\begin{proposition}[Regularized reproducing-kernel formula]
\label{prop:regularized-reproducing-kernel}
Fix $d\geq2$, $n\geq0$, and $x,y\in\BB^d$. Then the holomorphic extension of $\mathcal P_n^\mu(x,y)$ to the half-plane $\operatorname{Re} \mu>-1/2$ is given by
\begin{equation}
\label{eq:regularized-reproducing-kernel}
    \mathcal P_n^\mu(x,y)
    =2\int_{(0,1]}
    \left(\Theta_{n,E}^\mu(u;x,y)-\Theta_{n,E}^\mu(1;x,y)\right)
    \dint\Pi_{\mu-1/2}(u)
    +\Theta_{n,E}^\mu(1;x,y).
\end{equation}
In particular, if $-1/2<\mu<0$, then
\begin{equation}
\label{eq:subcritical-reproducing-kernel}
    \mathcal P_n^\mu(x,y)
    =-2\int_0^1\Pi_{\mu-1/2}(u) \, \partial_u\Theta_{n,E}^\mu(u;x,y) \, \dint u
    +\Theta_{n,E}^\mu(1;x,y).
\end{equation}
\end{proposition}

\begin{proof}
For real $\mu\geq0$, multiplying \eqref{eq:classical-reproducing-kernel} by $W_\mu(\BB^d)\lambda_\mu/(n+\lambda_\mu)$, we obtain
\[
    \mathcal P_n^\mu(x,y)=\int_{[-1,1]}\Theta_n^\mu(u;x,y) \, \dint\Pi_{\mu-1/2}(u).
\]
Since $\dint\Pi_{\mu-1/2}$ is even, only the even part of the integrand contributes, so that
\[
    \mathcal P_n^\mu(x,y)=2\int_{(0,1]}\Theta_{n,E}^\mu(u;x,y) \, \dint\Pi_{\mu-1/2}(u).
\]
Subtracting and adding the endpoint value $\Theta_{n,E}^\mu(1;x,y)$ gives \eqref{eq:regularized-reproducing-kernel} in the classical range.

It remains to extend the identity. By Proposition~\ref{prop:appendix-holomorphy}\textnormal{(i)} from Appendix~\ref{app:holomorphy}, the left-hand side of \eqref{eq:regularized-reproducing-kernel} admits a holomorphic extension to the half-plane $\operatorname{Re} \mu>-1/2$. Indeed, we use the identity \eqref{eq:renormalized-reproducing-kernel} and observe that the scalar factor on the right-hand side of \eqref{eq:renormalized-reproducing-kernel} is holomorphic in this half-plane, since $\lambda_\mu$ and $n+\lambda_\mu$ do not vanish there, while
\[
    W_\mu(\BB^d)=\pi^{d/2}\frac{\Gamma(\mu+1/2)}{\Gamma(\mu+(d+1)/2)}.
\]
Next, the endpoint term on the right-hand side of \eqref{eq:regularized-reproducing-kernel} is holomorphic because $C_n^{\lambda_\mu}$ is holomorphic in the present range of $\lambda_\mu$. Moreover, the function
\[
    \varphi_\mu(u)=\Theta_{n,E}^\mu(u;x,y)-\Theta_{n,E}^\mu(1;x,y)
\]
satisfies $\abs{\varphi_\mu(u)}\leq C_K(1-u)$ locally uniformly for $\mu$ in compact subsets $K$ of the half-plane. This follows from the mean value theorem, since $\partial_u\Theta_{n,E}^\mu$ is locally uniformly bounded. Proposition~\ref{prop:one-variable-Pi} gives holomorphy of the regularized integral, and the identity theorem proves \eqref{eq:regularized-reproducing-kernel}.

For real $-1/2<\mu<0$, formula \eqref{eq:subcritical-reproducing-kernel} follows from \eqref{eq:regularized-reproducing-kernel} by a standard integration-by-parts argument; see the proof of Proposition~\ref{prop:regularized-ball-reduction} below.
\end{proof}

\subsection{Analytic continuation of the ball reduction formula}
\label{sec:continued-reduction}

In this subsection, we record the regularized form of the ball heat-kernel reduction formula. For complex values of $\mu$, it is understood as an identity for the holomorphic extension provided by Proposition~\ref{prop:appendix-holomorphy}\textnormal{(ii)}. Since there is only one auxiliary variable, the formula has only one endpoint subtraction.

For $t>0$ and $x,y\in\BB^d$, define
\begin{equation}
\label{eq:G-mu-u}
    G_t^\mu(u;x,y) \coloneqq G_t^{\nu_\mu}(A+Bu,1),\qquad -1\leq u\leq1,
\end{equation}
and its even part
\begin{equation}
\label{eq:G-even}
    G_{t,E}^\mu(u;x,y) \coloneqq \frac12\left(G_t^\mu(u;x,y)+G_t^\mu(-u;x,y)\right),
    \qquad 0\leq u\leq1.
\end{equation}

\begin{proposition}[Regularized ball reduction formula]
\label{prop:regularized-ball-reduction}
Fix $d\geq2$, $t>0$, and $x,y\in\BB^d$. Then the holomorphic extension of $h_t^\mu(x,y)$ to the half-plane $\operatorname{Re}\mu>-1/2$ is given by
\begin{equation}
\label{eq:regularized-ball-reduction}
    h_t^\mu(x,y)
    =2c_{\mu,d}\int_{(0,1]}
    \left(G_{t,E}^\mu(u;x,y)-G_{t,E}^\mu(1;x,y)\right)
    \dint\Pi_{\mu-1/2}(u)
    +c_{\mu,d}G_{t,E}^\mu(1;x,y).
\end{equation}
In particular, if $-1/2<\mu<0$, then
\begin{equation}
\label{eq:subcritical-reduction}
    h_t^\mu(x,y)
    =-2c_{\mu,d}\int_0^1
    \Pi_{\mu-1/2}(u) \, \partial_uG_{t,E}^\mu(u;x,y) \, \dint u
    +c_{\mu,d}G_{t,E}^\mu(1;x,y).
\end{equation}
\end{proposition}

\begin{proof}
For real $\mu\geq0$, by \eqref{eq:classical-ball-reduction}, \eqref{eq:G-mu-u}, \eqref{eq:G-even}, and the symmetry of $\dint\Pi_{\mu-1/2}$, we obtain
\[
    h_t^\mu(x,y)=2c_{\mu,d}\int_{(0,1]}G_{t,E}^\mu(u;x,y) \, \dint\Pi_{\mu-1/2}(u).
\]
Subtracting and adding the boundary value $G_{t,E}^\mu(1;x,y)$ gives \eqref{eq:regularized-ball-reduction} in the classical range.

It remains to extend the identity. By Proposition~\ref{prop:appendix-holomorphy}\textnormal{(ii)} from Appendix~\ref{app:holomorphy}, the left-hand side of \eqref{eq:regularized-ball-reduction} admits a holomorphic extension to the half-plane $\operatorname{Re}\mu>-1/2$. 
Next, the endpoint term on the right-hand side of \eqref{eq:regularized-ball-reduction} is holomorphic there because
$c_{\mu,d}$ is holomorphic and the Jacobi heat series is locally
uniformly convergent and holomorphic in its parameters due to the fact that
$\operatorname{Re} \nu_\mu \geq \operatorname{Re} \mu >-1/2$. Moreover, the function
\[
    \varphi_\mu(u)=G_{t,E}^\mu(u;x,y)-G_{t,E}^\mu(1;x,y)
\]
satisfies
$
    \abs{\varphi_\mu(u)}\leq C_K(1-u),
$
locally uniformly for $\mu$ in compact subsets $K$ of the half-plane. This follows from the mean value theorem and the local uniform boundedness of the derivative in the auxiliary variable $u$.
By Proposition~\ref{prop:one-variable-Pi}, the regularized integral is holomorphic in $\mu$. The identity theorem gives \eqref{eq:regularized-ball-reduction} for all $\mu$ such that $\operatorname{Re}\mu>-1/2$, with the left-hand side understood as the holomorphic extension when $\mu$ is not real.

Finally, let $-1/2<\mu<0$ and put $\alpha \coloneqq \mu-1/2\in(-1,-1/2)$. Set
\[
    F(u) \coloneqq G_{t,E}^\mu(u;x,y)-G_{t,E}^\mu(1;x,y).
\]
Then $F\in C^1([0,1])$ and $F(1)=0$. Integrating by parts on $[0,1-\varepsilon]$ gives
\[
    \int_0^{1-\varepsilon}F(u) \, \dint\Pi_\alpha(u)
    =F(1-\varepsilon)\Pi_\alpha(1-\varepsilon)-F(0)\Pi_\alpha(0)
    -\int_0^{1-\varepsilon}\Pi_\alpha(u)F'(u) \, \dint u.
\]
The lower boundary term vanishes because $\Pi_\alpha(0)=0$. The upper boundary term tends to zero because $F(u)=O(1-u)$ and $\abs{\Pi_\alpha(u)}\lesssim(1-u)^{\alpha+1/2}$ near $u=1$ by \eqref{eq:Pi-size}. Since $\alpha+1/2>-1$, the last integral converges absolutely. Passing to the limit gives
\[
    \int_{(0,1]}F(u) \, \dint\Pi_\alpha(u)
    =-\int_0^1\Pi_\alpha(u)F'(u) \, \dint u,
\]
and this yields \eqref{eq:subcritical-reduction}.
\end{proof}

\begin{remark}
From now on, we mainly use formula \eqref{eq:subcritical-reduction}, where the singular $\dint\Pi_{\mu-1/2}$ pairing is replaced by the sum of an endpoint term and an ordinary derivative correction integral.
\end{remark}

\subsection{Sharp estimates on the ball}
\label{sec:ball-proof}

We now prove Theorem~\ref{thm:main-ball}. Throughout this subsection, we fix $-1/2 < \mu < 0$ and $\nu \coloneqq \nu_\mu=\mu+\frac{d-2}{2}$. Given $x,y \in \BB^d$, we put
\begin{equation*}
    d_{\BB} \coloneqq d_{\BB}(x,y), \qquad q \coloneqq \frac{B}{\pi-d_{\BB}}=\frac{\delta(x)\delta(y)}{\pi-d_{\BB}(x,y)}.
\end{equation*}
We may treat the points $x, y$ as fixed, while requiring all constants to be independent of them. 

By \eqref{eq:subcritical-reduction}, we split
\begin{equation}
\label{eq:D-E-split}
    h_t^\mu(x,y)=D_t(x,y)+E_t(x,y),
\end{equation}
where
\begin{equation}
\label{eq:D-term}
    D_t(x,y) \coloneqq -2c_{\mu,d}\int_0^1
    \Pi_{\mu-1/2}(u) \, \partial_uG_{t,E}^\mu(u;x,y) \, \dint u
\end{equation}
and
\begin{equation}
\label{eq:E-term}
    E_t(x,y) \coloneqq c_{\mu,d}G_{t,E}^\mu(1;x,y).
\end{equation}

\subsubsection{Auxiliary estimates}

The derivative identity and the sharp one-point estimate for the ultraspherical heat kernel imply the following bound.

\begin{lemma}[Derivative comparison]
\label{lem:derivative-comparison}
Fix $T>0$. Then one has
\begin{equation}
\label{eq:derivative-comparison}
    \abs{\partial_zG_t^\nu(z,1)}
    \lesssim
    t^{-\nu-2}[t+\pi-\arccos z]^{-\nu-3/2}
    \exp\left(-\frac{\arccos(z)^2}{4t}\right),
\end{equation}
uniformly in $z\in[-1,1]$ and $0<t\leq T$. Consequently, denoting
\[
    \Phi(u)=\arccos(A+Bu)
\]
for $u \in [0,1]$, one has
\begin{equation}
\label{eq:u-derivative-comparison}
    \abs{\partial_uG_{t,E}^\mu(u;x,y)}
    \lesssim
    B t^{-\nu-2}[t+\pi-\Phi(u)]^{-\nu-3/2}
    \exp\left(-\frac{\Phi(u)^2}{4t}\right).
\end{equation}
\end{lemma}

\begin{proof}
    The estimates \eqref{eq:derivative-comparison} and \eqref{eq:u-derivative-comparison} follow readily from \eqref{eq:derivative-ultra} and \eqref{eq:one-point-ultra}.
\end{proof}

The next lemma separates the two branches generated by the even part.

\begin{lemma}[Geometric separation]
\label{lem:geometric-separation}
Let
\[
    \Phi_-(u)=\arccos(A-Bu),\qquad 0\leq u\leq1.
\]
Then one has
\begin{equation}
\label{eq:geometric-separation}
    \Phi_-(u)^2-d_{\BB}^2\gtrsim \frac{B}{\pi-d_{\BB}},
\end{equation}
where the quotient $B / (\pi-d_{\BB})$ is understood as $0$ when $d_\BB = \pi$. Consequently, one has
\begin{equation}
\label{eq:gaussian-separation}
    \exp\left(-\frac{\Phi_-(u)^2}{4t}\right)
    \leq
    \exp\left(-\frac{d_{\BB}^2}{4t}\right)
    \exp\left(-\frac{cB}{t(\pi-d_{\BB})}\right)
\end{equation}
with an absolute constant $c>0$.
\end{lemma}

\begin{proof}
If $d_{\BB} = 0$ or $d_{\BB} = \pi$, then \eqref{eq:geometric-separation} readily holds. Otherwise, put $x_0=A+B=\cos d_{\BB}$ and $\eta=x_0-(A-Bu)=B(1+u)$. Let $f(z)=\arccos(z)^2$. By the mean value theorem,
\[
    \Phi_-(u)^2-d_{\BB}^2=f(x_0-\eta)-f(x_0)=-\eta f'(\xi)
\]
for some $\xi\in(x_0-\eta,x_0)$. Observe that the expression
\[
    -f'(z)=\frac{2\arccos z}{\sqrt{1-z^2}}
\]
is decreasing in $z\in(-1,1)$. Indeed, substituting $\arccos z \mapsto \theta$, we obtain
\[
\frac{2\arccos z}{\sqrt{1-z^2}} = \frac{2\theta}{\sin \theta}
\]
and the right-hand side is clearly increasing in $\theta \in (0,\pi)$. Therefore, we get
\[
    \Phi_-(u)^2-d_{\BB}^2 = -\eta f'(\xi) \geq
    -\eta f'(x_0) = \eta\frac{2d_{\BB}}{\sin d_{\BB}}.
\]
Since $\eta\geq B$, the elementary inequality
\[
    \frac{2s}{\sin s}\gtrsim \frac1{\pi-s},\qquad 0<s<\pi,
\]
applied with $s=d_{\BB}$ gives \eqref{eq:geometric-separation}. The Gaussian estimate \eqref{eq:gaussian-separation} follows immediately.
\end{proof}

\begin{lemma}[Size bound for $q$]
\label{lem:q-size-bound}
One has
\begin{equation}
\label{eq:q-size-bound}
    q\lesssim \pi-d_{\BB}.
\end{equation}
In particular, if $q\gtrsim t$, then $t+\pi-d_{\BB}\simeq \pi-d_{\BB}$.
\end{lemma}

\begin{proof}
Since $A-B\geq -1$, we have $1+A+B\geq 2B$. Consequently,
\[
    B\leq \frac{1+A+B}{2} = \frac{1+\cos d_{\BB}}{2}=\cos^2\frac{d_{\BB}}2=\sin^2\frac{\pi-d_{\BB}}{2}\lesssim (\pi-d_{\BB})^2
\]
and dividing by $\pi-d_{\BB}$ gives \eqref{eq:q-size-bound}. The final assertion follows immediately.
\end{proof}

\begin{lemma}[Branch dominance]
\label{lem:branch-dominance}
There exists a large constant $M>1$ such that if
\[
    q=\frac{B}{\pi-d_{\BB}}\geq Mt
\]
and
\[
    I_* = [1- t/q,1 ] \subset (0,1],
\]
then, with $\Phi(u)=\arccos(A+Bu)$, for every $u\in I_*$ one has
\begin{equation}
\label{eq:forward-square-localization}
    0\leq \Phi(u)^2-d_{\BB}^2\lesssim t
\end{equation}
and
\begin{equation}
\label{eq:forward-antipodal-comparison}
    t+\pi-\Phi(u)\simeq t+\pi-d_{\BB}.
\end{equation}
Moreover,
\begin{equation}
\label{eq:branch-dominance}
    G_t^{\nu+1}(A-Bu,1)
    \leq \frac12 G_t^{\nu+1}(A+Bu,1).
\end{equation}
\end{lemma}

\begin{proof}
Put $\Delta=\pi-d_{\BB}$. Since $q\geq Mt>0$, we have
$d_{\BB}<\pi$ and $B>0$. For $u\in I_*$, we have
\begin{equation}
\label{eq:forward-cosine-defect}
    \cos d_{\BB} -\cos\Phi(u)=
    (A+B) - (A + Bu) = B(1-u) \leq B t/q = t\Delta.
\end{equation}
Assume that $0\leq d_\BB \leq\phi\leq\pi$ and observe that 
$\phi-d_\BB \lesssim\sin((\phi-d_\BB)/2)$. We also have 
\[
(\phi+d_\BB)(\pi-d_\BB)\lesssim\sin((\phi+d_\BB)/2),
\]
where the case $\phi+d_\BB \leq\pi$ is standard and the case $\phi+d_\BB >\pi$ follows by applying the bound
$\pi-d_\BB \leq 2\pi-(\phi+d_\BB) \lesssim \sin(\pi-(\phi+d_\BB)/2) = \sin((\phi+d_\BB)/2)$. This leads to the estimate
\[
    \phi^2-d_\BB^2\lesssim
    \frac{2 \sin((\phi-d_\BB)/2) \sin((\phi+d_\BB)/2)}{\pi-d_\BB} = \frac{\cos d_\BB-\cos\phi}{\pi-d_\BB},
\]
which applied with $\phi=\Phi(u) \in [d_\BB, \pi]$ gives
\eqref{eq:forward-square-localization}, as desired.

We next prove \eqref{eq:forward-antipodal-comparison}. Since $\Phi(u) \in [d_\BB, \pi]$, it suffices to consider the case $d_\BB > \pi/2$. We put
$\varepsilon = \pi-\Phi(u)$ and note that $\varepsilon \leq\Delta<\pi/2$. By the previous estimate, we have
\[
    \cos d_{\BB} -\cos\Phi(u) 
    \gtrsim (\pi - d_\BB)(\Phi(u)^2 - d_\BB^2) \simeq \Delta(\Delta-\varepsilon),
\]
which gives $0\leq\Delta-\varepsilon\lesssim t$ by \eqref{eq:forward-cosine-defect}. This completes the proof of \eqref{eq:forward-antipodal-comparison}. The latter argument remains valid when $d_\BB=0$, as no comparison between $\Phi(u)$ and $d_\BB$ is used in this case.

Let $\Phi_-(u)=\arccos(A-Bu)$. Lemma~\ref{lem:geometric-separation} and
\eqref{eq:forward-square-localization} imply, after adjusting $M$, that
\begin{equation}
\label{eq:branch-square-gap}
    \Phi_-(u)^2-\Phi(u)^2\gtrsim q.
\end{equation}
The same cosine comparison, together with
$\cos(\pi-\Phi_-(u))-\cos\Delta=B(1+u)\leq2q\Delta$, gives
\[
    \Delta\lesssim \pi-\Phi_-(u)+q.
\]
Indeed, if $\Delta\leq\pi/2$ then we use
$\cos(\pi-\Phi_-(u))-\cos\Delta
\gtrsim \Delta[\Delta-(\pi-\Phi_-(u))]$ much as before, while if $\Delta>\pi/2$, then we can consider only the case
$\pi-\Phi_-(u)<\Delta/2$ and conclude that the above difference of cosines exceeds some absolute constant, and so $q\gtrsim1$. Consequently,
\[
    \frac{t+\pi-\Phi(u)}{t+\pi-\Phi_-(u)}
    = 1 + \frac{\Phi_-(u) - \Phi(u)}{t+\pi-\Phi_-(u)}
    \leq 1 + \frac{\Delta}{t+\pi-\Phi_-(u)}
    \lesssim 1+\frac qt.
\]
Using the two-sided estimate \eqref{eq:one-point-ultra} with parameter
$\nu+1$ and \eqref{eq:branch-square-gap}, we obtain
\[
    \frac{G_t^{\nu+1}(A-Bu,1)}{G_t^{\nu+1}(A+Bu,1)}
    \lesssim
    \left(1+\frac qt\right)^{\nu+3/2}
    \exp\left(-c\frac qt\right).
\]
Choosing $M$ sufficiently large makes the above quotient at most $1/2$ and
proves \eqref{eq:branch-dominance}.
\end{proof}

\subsubsection{Upper bound}

By \eqref{eq:D-E-split}, it suffices to estimate \eqref{eq:D-term} and \eqref{eq:E-term}.

We first estimate the endpoint term \eqref{eq:E-term}. By monotonicity, we have
\[
    E_t(x,y)=\frac{c_{\mu,d}}2(G_t^\nu(A+B,1)+G_t^\nu(A-B,1)) \lesssim G_t^\nu(A+B,1).
\]
Moreover, \eqref{eq:one-point-ultra} gives
\[
    G_t^\nu(A+B,1)
    \lesssim
    t^{-\nu-1}[t+\pi-d_{\BB}]^{-\nu-1/2}
    \exp\left(-\frac{d_{\BB}^2}{4t}\right).
\]
We note that $-\nu-1=-d/2-\mu$ and $-\nu-1/2= - \lambda_\mu$ and conclude that the right-hand side above is dominated by the right-hand side of \eqref{eq:main-estimate} when $-1/2 < \mu < 0$.

For the derivative term \eqref{eq:D-term}, we first use Lemma~\ref{lem:derivative-comparison}. By \eqref{eq:u-derivative-comparison}, we obtain
\[
|D_t(x,y)| \lesssim B t^{-\nu-2} \int_0^1
    [t+\pi-\Phi(u)]^{-\nu-3/2}
    \exp\left(-\frac{\Phi(u)^2}{4t}\right) \, |\Pi_{\mu-1/2}(u)| \, \dint u.
\]
Using \eqref{eq:Pi-positive}, \eqref{eq:Pi-size}, and $|u| \lesssim 1$, we can further estimate
\[
|D_t(x,y)| \lesssim B t^{-\nu-2}
\int_{[0,1]}
    [t+\pi-\Phi(u)]^{-\nu-3/2}
    \exp\left(-\frac{\Phi(u)^2}{4t}\right) \, \dint\Pi_{\mu+1/2}(u).
\]
Applying Lemma~\ref{lem:integral-engine} with $\gamma=\mu+1/2$, $\xi=\nu +3/2$, and $D=4t \lesssim t$, we obtain
\[
|D_t(x,y)| \lesssim B t^{-\nu-2+\mu+1}
    [t+\pi-d_{\BB}]^{-\nu-3/2} (t+q)^{-\mu-1}
    \exp\left(-\frac{d_{\BB}^2}{4t}\right).
\]
Since $B=q(\pi-d_{\BB})\leq (t+q)(t+\pi-d_{\BB})$, this gives
\[
|D_t(x,y)| \lesssim t^{-\nu-2+\mu+1}
    [t+\pi-d_{\BB}]^{-\nu-1/2} (t+q)^{-\mu} 
    \exp\left(-\frac{d_{\BB}^2}{4t}\right).   
\]
We note that $-\nu-2+\mu+1 = -d/2$ and conclude that the right-hand side above is dominated by the right-hand side of \eqref{eq:main-estimate}, as desired. 

\subsubsection{Lower bound}

The lower bound is obtained by splitting the proof into two cases.

\paragraph{{\bf Case 1: $q\gtrsim t$.}}

Assume that $q\geq Mt$, with $M>1$ as in Lemma~\ref{lem:branch-dominance}. By \eqref{eq:derivative-ultra}, the derivative term \eqref{eq:D-term} can be written in the form
\begin{align}
\label{eq:D-positive-form}
    D_t(x,y)
    &=C_{\mu,d,t}B\int_0^1\abs{\Pi_{\mu-1/2}(u)}
    \big(G_t^{\nu+1}(A+Bu,1)-G_t^{\nu+1}(A-Bu,1)\big) \, \dint u
\end{align}
for $0<t\leq T$, with some non-negative constant $C_{\mu,d,t}\simeq1$.

Consider the interval
\[
    I_* = [1- t/q,1]
\]
from Lemma~\ref{lem:branch-dominance}. For all $u \in [0,1]$, the integrand in \eqref{eq:D-positive-form} is non-negative. By Lemma~\ref{lem:branch-dominance}, if $u \in I_*$, then the lower bound for the difference in \eqref{eq:D-positive-form} is $\frac12 G_t^{\nu+1}(A+Bu,1)$. This gives
\begin{align*}
    D_t(x,y)
    &\gtrsim B\int_{I_*}\abs{\Pi_{\mu-1/2}(u)} \, G_t^{\nu+1}(A+Bu,1) \, \dint u.
\end{align*}
By the two-sided estimate \eqref{eq:one-point-ultra} with parameter $\nu+1$, followed by \eqref{eq:forward-square-localization} and
\eqref{eq:forward-antipodal-comparison}, we have
\begin{align}
\label{eq:D-local-lower}
    D_t(x,y)
    &\gtrsim
    B t^{-\nu-2}[t+\pi-d_{\BB}]^{-\nu-3/2}
    \exp\left(-\frac{d_{\BB}^2}{4t}\right)
    \int_{I_*}\abs{\Pi_{\mu-1/2}(u)} \, \dint u.
\end{align}
Since $\mu-1/2\in(-1,-1/2)$, we have
\[
    \abs{\Pi_{\mu-1/2}(u)}\simeq(1-u)^\mu, \qquad u \in I_*.
\]
The length of $I_*$ equals $t/q$. Thus,
\begin{equation}
\label{eq:Istar-weight}
    \int_{I_*}\abs{\Pi_{\mu-1/2}(u)} \, \dint u
    \simeq \left(\frac{t}{q}\right)^{\mu+1}.
\end{equation}
Substituting \eqref{eq:Istar-weight} into \eqref{eq:D-local-lower}, using the estimate $B/q = \pi - d_\BB \simeq t+\pi - d_\BB$ for $q\gtrsim t$ from Lemma~\ref{lem:q-size-bound}, and simplifying, we conclude that
\[
    D_t(x,y)
    \gtrsim
    [t+\pi-d_{\BB}]^{-\lambda_\mu}
    q^{-\mu}
    t^{-d/2}
    \exp\left(-\frac{d_{\BB}^2}{4t}\right).
\]
For $q\gtrsim t$ we have $(t+q)^{-\mu}\simeq q^{-\mu}$. By $E_t(x,y) \geq 0$ the desired lower bound holds in this case.

\paragraph{{\bf Case 2: $q\lesssim t$.}}

As in the previous case, we observe that $D_t(x,y) \geq 0$. In this case, we have
\[
    (t+q)^{-\mu}\simeq t^{-\mu}
\]
and the endpoint term \eqref{eq:E-term} provides the required lower bound. Indeed, we have
\[
    E_t(x,y)\geq cG_t^\nu(A+B,1)=cG_t^\nu(\cos d_{\BB},1),
\]
and by \eqref{eq:one-point-ultra} this gives
\[
    E_t(x,y)
    \gtrsim
    t^{-\nu-1}[t+\pi-d_{\BB}]^{-\nu-1/2}
    \exp\left(-\frac{d_{\BB}^2}{4t}\right).
\]
Since $-\nu-1=-d/2-\mu$, this is exactly the desired bound in the case $q\lesssim t$.

Combining the two cases proves the lower estimate in \eqref{eq:main-estimate}. The upper and lower estimates together prove Theorem~\ref{thm:main-ball} for $0<t\leq T$.

The estimate $h_t^\mu(x,y)\simeq1$ for $t\geq T$ follows from the spectral expansion \eqref{eq:ball-heat-series}, the positivity of the semigroup, and the fact that the eigenvalue $0$ is simple and corresponds to constants. 

This completes the proof of Theorem~\ref{thm:main-ball}.

\section{The simplex case}
\label{sec:simplex-case}

\subsection{A regularized formula for the simplex reproducing kernel}
\label{sec:regularized-simplex-reproducing-kernel}

We now repeat the idea from
Subsection~\ref{sec:regularized-reproducing-kernel}, one coordinate at a
time. For $\operatorname{Re}\alpha>-1$ and $f\in C^1([-1,1])$, define
\begin{equation}
\label{eq:regularization-operator}
    \mathfrak R_\alpha f
    \coloneqq
    \frac{f(-1)+f(1)}2
    -\int_{-1}^1\Pi_\alpha(u)f'(u) \, \dint u,
\end{equation}
with the integral omitted for $\alpha = - \frac12$. When
$\operatorname{Re}\alpha>-1/2$, integration by parts gives
\[
    \mathfrak R_\alpha f
    =\int_{[-1,1]}f(u) \, \dint\Pi_\alpha(u),
\]
since $\Pi_\alpha(1) = \frac{1}{2}$ and $\Pi_\alpha(-1) = -\frac{1}{2}$. This formula also justifies our convention for $\alpha = - \frac12$.

Proposition~\ref{prop:one-variable-Pi}, applied after separating the two
halves of the interval, shows that \eqref{eq:regularization-operator} is
the holomorphic extension of this pairing to $\operatorname{Re}\alpha>-1$.
For a function of several variables, $\mathfrak R_{\alpha_j}^{(j)}$
means that this operation is performed in the $j$th variable. The operations
commute.

Put
\begin{equation}
\label{eq:simplex-normalizing-constant}
    c_{\kappa,d}^{\VV}
    \coloneqq
    \frac{\sqrt\pi\,\Gamma(\rho_\kappa+2)}
    {\prod_{j=1}^{d+1}\Gamma(\kappa_j+1/2)}
\end{equation}
and
\begin{equation}
\label{eq:simplex-Psi}
    \Psi_n^\kappa(u;x,y)
    \coloneqq
    \mathfrak D_{2n}^{\rho_\kappa}
    P_{2n}^{\rho_\kappa}\big(z(u;x,y)\big).
\end{equation}

\begin{proposition}[Regularized simplex reproducing formula]
\label{prop:regularized-simplex-reproducing-kernel}
Let $d\geq2$, $\kappa\in(-1/2,\infty)^{d+1}$, $n\geq0$, and
$x,y\in\VV^d$. Then
\begin{equation}
\label{eq:regularized-simplex-reproducing-kernel}
    P_n(U_\kappa,x,y)
    =c_{\kappa,d}^{\VV}
    \left(\prod_{j=1}^{d+1}
    \mathfrak R_{\kappa_j-1/2}^{(j)}\right)
    \Psi_n^\kappa(\,\cdot\,;x,y).
\end{equation}
More explicitly, let $N=N(\kappa)=\{j:\kappa_j<0\}$ and $P=\{1,\ldots,d+1\}\setminus N$. 
Denoting $\partial_S=\prod_{j\in S}\partial_{u_j}$ for $S \subset N$, we can expand the right-hand side of
\eqref{eq:regularized-simplex-reproducing-kernel} as
\begin{align*}
    c_{\kappa,d}^{\VV}
    \sum_{S\subset N}(-1)^{\abs{S}}
    \int_{[-1,1]^{d+1}}
    \partial_S\Psi_n^\kappa(u;x,y)
    \prod_{j\in S}\Pi_{\kappa_j-1/2}(u_j) \, \dint u_j 
    \prod_{j\in N\setminus S}\dint\Pi_{-1/2}(u_j)
    \prod_{j\in P}\dint\Pi_{\kappa_j-1/2}(u_j).
\end{align*}
\end{proposition}

\begin{proof}
In the classical range $\kappa\in[0,\infty)^{d+1}$, the original simplex product formula,
followed by the quadratic Jacobi transformation, gives
\[
    P_n(U_\kappa,x,y)
    =c_{\kappa,d}^{\VV}
    \int_{[-1,1]^{d+1}}
    \Psi_n^\kappa(u;x,y)
    \prod_{j=1}^{d+1}\dint\Pi_{\kappa_j-1/2}(u_j).
\]
This is the identity underlying the simplex reduction in \cite[Section~6]{NSS19}. In the extended range $\kappa\in(-1/2,\infty)^{d+1}$, we want to apply \eqref{eq:regularization-operator} in the coordinates $j$ with $\kappa_j < 0$.

 The factors in
\eqref{eq:simplex-normalizing-constant} and
\eqref{eq:simplex-Psi} are holomorphic throughout the complex domain
$\mathcal T_{\VV}$ of Appendix~\ref{app:holomorphy}; the use of even
indices is essential here because
$\mathfrak D_{2n}^\rho$ extends holomorphically to $\operatorname{Re}\rho>-3/2$.
Applying \eqref{eq:regularization-operator} successively in the
$d+1$ variables gives a jointly holomorphic right-hand side. The simplex
analogue of Proposition~\ref{prop:appendix-holomorphy}, proved in
Proposition~\ref{prop:appendix-simplex-holomorphy}, ensures that the left-hand side is holomorphic in each variable separately. Continuing one parameter
at a time and using the identity theorem proves
\eqref{eq:regularized-simplex-reproducing-kernel}. 
\end{proof}

\subsection{Analytic continuation of the simplex reduction formula}
\label{sec:continued-simplex-reduction}

The heat-kernel counterpart of
Proposition~\ref{prop:regularized-simplex-reproducing-kernel} has exactly
the same form.

\begin{proposition}[Regularized simplex reduction formula]
\label{prop:regularized-simplex-reduction}
Let $d\geq2$, $\kappa\in(-1/2,\infty)^{d+1}$, $t>0$, and
$x,y\in\VV^d$. Then
\begin{equation}
\label{eq:regularized-simplex-reduction}
    H_t^\kappa(x,y)
    =c_{\kappa,d}^{\VV}
    \left(\prod_{j=1}^{d+1}
    \mathfrak R_{\kappa_j-1/2}^{(j)}\right)
    \mathcal H_{t/4}^{\rho_\kappa}
    \big(z(\,\cdot\,;x,y)\big).
\end{equation}
Equivalently, with $N$ and $P$ as above, we can write
\begin{equation}
\label{eq:simplex-subset-decomposition}
    H_t^\kappa(x,y)=\sum_{S\subset N}D_{t,S}^\kappa(x,y),
\end{equation}
where each $D_{t,S}^\kappa(x,y)$ is given by 
\begin{align*}
    c_{\kappa,d}^{\VV}(-1)^{\abs S}
    \int_{[-1,1]^{d+1}}
    \partial_S\!\left[
    \mathcal H_{t/4}^{\rho_\kappa} \big(z(u;x,y)\big)\right]
    \prod_{j\in S}\Pi_{\kappa_j-1/2}(u_j) \,\dint u_j 
    \prod_{j\in N\setminus S}\dint\Pi_{-1/2}(u_j)
    \prod_{j\in P}\dint\Pi_{\kappa_j-1/2}(u_j).
\end{align*}
Here $D_{t,\varnothing}^\kappa$ is the endpoint contribution, while the other terms are mixed derivative corrections.
\end{proposition}

\begin{proof}
We multiply \eqref{eq:regularized-simplex-reproducing-kernel} by
$\me^{-tn(n+\lambda_\kappa)}$ and sum in $n$. Since
\[
    tn(n+\lambda_\kappa)
    =\frac t4\,2n(2n+2\rho_\kappa+1),
\]
the resulting series is
\eqref{eq:auxiliary-H-rho}. Local uniform convergence follows from Jacobi-polynomial and gamma-ratio bounds analogous to those used in Appendix~\ref{app:holomorphy}. This proves
\eqref{eq:regularized-simplex-reduction}, and expanding the regularization
operators gives \eqref{eq:simplex-subset-decomposition}.
\end{proof}

\begin{remark}
When $N$ has one element, \eqref{eq:simplex-subset-decomposition} is the
same endpoint-plus-derivative decomposition as in
\eqref{eq:subcritical-reduction}. For general $N$, its $2^{\abs N}$ terms
record independently whether each singular coordinate contributes its
endpoint value or its derivative correction.
\end{remark}

\subsection{Sharp estimates on the simplex}
\label{sec:simplex}

We now prove Theorem~\ref{thm:main-simplex}. Throughout this subsection, we fix $\kappa \in (-1/2, \infty)^{d+1}$. Given $x,y \in \VV^d$, we put
\[
    d_\VV=d_{\VV}(x,y),
    \qquad
    \rho=\rho_\kappa.
\]
We may treat the points $x, y$ as fixed, while requiring all constants to be independent of them.

We observe that the absence of an antipodal factor in
\eqref{eq:main-simplex-estimate} reflects the elementary fact that the positive cube $[0,1]^{d+1}$ only produces angles in $[0,\pi/2]$.

\subsubsection{Auxiliary estimates}

We first record the iterated form of
Lemma~\ref{lem:integral-engine} used below. 

\begin{lemma}[Iterated simplex integral estimate]
\label{lem:iterated-simplex-integral}
Let $J \subset \{1,\dots,d+1\}$ be a finite index set. Fix $\gamma_j\geq-1/2$, $b_j\geq0$, and
$A\geq0$ such that $A+\sum_{j\in J}b_j\leq1$. Then
\begin{align}
\begin{split}
    \label{eq:iterated-simplex-integral}
    &\int_{[0,1]^{J}}
    \exp\left(-\frac{\arccos(A+\sum_{j\in J}b_ju_j)^2}{t}\right)
    \prod_{j\in J}\dint\Pi_{\gamma_j}(u_j) \\
    &\quad\simeq
    t^{\sum_{j\in J}(\gamma_j+1/2)}
    \prod_{j\in J}(b_j + t)^{-\gamma_j-1/2}
    \exp\left(-\frac{\arccos(A+\sum_{j\in J}b_j)^2}{t}\right),
\end{split}
\end{align}
uniformly for $0<t\leq T$. Here $[0,1]^J$ stands for the cube $[0,1]^{|J|}$ with coordinates indexed by $J$.
\end{lemma}

\begin{proof}
We apply Lemma~\ref{lem:integral-engine} successively in the variables from
$J$, with $\xi=0$. At every step the terminal angle $\arccos(A+B)$ belongs to
$[0,\pi/2]$, so the factors $\pi- \arccos(A+B)$ can be ignored. This gives
\eqref{eq:iterated-simplex-integral}.
\end{proof}

The next observation is the multi-variable version of
Lemma~\ref{lem:branch-dominance}. Recall from \eqref{eq:simplex-a-z} that $z(u;x,y) = \sum_{j=1}^{d+1}a_ju_j$, where $a_j = \sqrt{x_jy_j}$. If $0\leq z_+\leq1$ is a sum of non-negative terms $a_j u_j$ and $z_-$ is
obtained from $z_+$ by changing the signs of some terms whose total size is
$b$, then
\begin{equation}
\label{eq:simplex-gaussian-separation}
    \arccos(z_-)^2-\arccos(z_+)^2\gtrsim b.
\end{equation}
Indeed, the derivative of $s\mapsto\arccos(s)^2$ has absolute value
$2\arccos(s)/\sqrt{1-s^2}\geq2$. 

Consequently, for $M > 0$ large enough, sign changes in terms whose total size is at least $Mt$ produce an additional factor $\me^{-cM}$ in the
Gaussian. As in the proof of \cite[Theorem~6.1]{NSS19}, such exponential
loss absorbs all polynomial factors coming from the one-point estimate.

Finally, every term in \eqref{eq:simplex-subset-decomposition} is
non-negative. To see this, we reflect all integrations to the positive
cube. Since $\Pi_{\kappa_j-1/2}(u)<0$ for $u\in(0,1)$ and $j\in S$,
the signs outside the integral cancel. The remaining signed sum under the integral can be written as
\[
\sum_{R \subset S} (-1)^{|R|} \partial_S\!\left[
    \mathcal H_{t/4}^{\rho_\kappa} \big(z(u^{(R)};x,y)\big)\right]
\]
for $u^{(R)}_j \coloneqq -u_j$ if $j \in R$ and $u^{(R)}_j \coloneqq u_j$ otherwise. By the fundamental
theorem of calculus, it can be rewritten as an average of the $2\abs S$th derivative.
Formula \eqref{eq:auxiliary-H-derivatives} expresses this derivative as
a positive multiple of the even part of
$G_{t/4}^{\rho+2\abs S}(\, \cdot \,,1)$. Hence, for each $S \subset N$, we have
\begin{equation}
\label{eq:simplex-D-positive}
    D_{t,S}^\kappa(x,y)\geq0.
\end{equation}

\subsubsection{Upper bound}

Let $S\subset N$ and $r=\abs S$. Reflecting to the positive cube,
using \eqref{eq:simplex-gaussian-separation}, and applying
\eqref{eq:auxiliary-H-estimate} when $r=0$ and
\eqref{eq:auxiliary-H-derivatives} when $r\geq1$, we can control $D_{t,S}^\kappa(x,y)$ by
\begin{align*}
    t^{-\rho-r-1}
    \int_{[0,1]^{P\cup S}}
    \exp\left(-\frac{\arccos(A+\sum_{j\in P\cup S}a_ju_j)^2}{t}\right)
    \prod_{j\in P}\dint\Pi_{\kappa_j-1/2}(u_j)
    \prod_{j\in S} a_j \abs{\Pi_{\kappa_j-1/2}(u_j)} \, \dint u_j,
\end{align*}
where $A=\sum_{j\in N\setminus S}a_j$. We replace $\abs{\Pi_{\kappa_j-1/2}(u_j)} \, \dint u_j$ by $\dint\Pi_{\kappa_j+1/2}(u_j)$.
Lemma~\ref{lem:iterated-simplex-integral} now gives
\begin{align*}
    D_{t,S}^\kappa(x,y)
    \lesssim{}&
    t^{-d/2-\sum_{j\in N\setminus S}\kappa_j}
    \prod_{j\in P}(t+a_j)^{-\kappa_j}
    \prod_{j\in S}
    \frac{a_j}{t+a_j}(t+a_j)^{-\kappa_j}
    \exp\left(-\frac{d_\VV^2}{t}\right).
\end{align*}
Since $a_j/(t+a_j)\leq1$ and, for $j\in N\setminus S$, we have $t^{-\kappa_j}\leq(t+a_j)^{-\kappa_j}$,
the right-hand side above is controlled by the right-hand side of
\eqref{eq:main-simplex-estimate}. Summing over $S\subset N$ proves the
upper bound.

\subsubsection{Lower bound}

Choose a large constant $M>1$ to be determined and set
\begin{equation*}
    L=\{j\in N:a_j\geq Mt\},
    \qquad R=N\setminus L.
\end{equation*}
By Proposition~\ref{prop:regularized-simplex-reduction} and \eqref{eq:simplex-D-positive}, it suffices to show that $D_{t,L}^\kappa$ has the desired lower bound. 

Suppose first that $L=\varnothing$. In the endpoint term
$D_{t,\varnothing}^\kappa$, choose the positive endpoint in every
coordinate from $N$ and restrict the remaining integrations to the
positive cube. By positivity, \eqref{eq:auxiliary-H-estimate}, and
Lemma~\ref{lem:iterated-simplex-integral}, we have
\begin{align*}
    D_{t,\varnothing}^\kappa(x,y)
    &\gtrsim
    t^{-d/2-\sum_{j\in N}\kappa_j}
    \prod_{j\in P}(t+a_j)^{-\kappa_j}
    \exp\left(-\frac{d_\VV^2}{t}\right).
\end{align*}
Since $a_j<Mt$ for $j\in N$, this is the required lower bound.

Assume now that $L\neq\varnothing$ and put $r=\abs L$. For $j\in L$,
let
\[
    I_j=[1- t / a_j ,1].
\]
In the symmetrized mixed
derivative occurring in $D_{t,L}^\kappa$, the branch with all positive
signs dominates the sum of the other branches on
$\prod_{j\in L}I_j$. Indeed,
\eqref{eq:simplex-gaussian-separation} gives a loss $\me^{-cM}$ for
each other branch, so it suffices to take $M>0$ large enough. 

Next, we choose the positive
endpoints in the coordinates from $R$ and the positive cube in the
coordinates from $P$. For the branch with all positive signs, placing the variables from $P$ at their
positive endpoints, we put
\[
    z_L(\boldsymbol u_L)
    =\sum_{j\in R}a_j+\sum_{j\in L}a_ju_j+\sum_{j\in P}a_j.
\]
Then, on $\prod_{j\in L}I_j$, we obtain
\[
    0\leq \cos d_\VV -z_L(\boldsymbol u_L)
    =\sum_{j\in L}a_j(1-u_j)
    \leq r t.
\]
Since the absolute value of the derivative of
$s\mapsto\arccos(s)^2$ is bounded on $[0,1]$, this gives
\[
    0\leq \arccos(z_L(\boldsymbol u_L))^2-d_\VV^2\lesssim t.
\]
Thus, the corresponding exponential term is uniformly comparable to
$\exp(-d_\VV^2/t)$ on $\prod_{j\in L}I_j$. 

Now, by
\eqref{eq:auxiliary-H-derivatives}, the discussion leading to \eqref{eq:simplex-D-positive}, and the one-point estimate, we have
\begin{align*}
    D_{t,L}^\kappa(x,y)
    \gtrsim{}
    t^{-\rho-r-1}
    \prod_{j\in L}
    a_j \int_{I_j}\abs{\Pi_{\kappa_j-1/2}(u_j)} \, \dint u_j
    \prod_{j\in P} t^{\kappa_j} (t+a_j)^{-\kappa_j}
    \exp\left(-\frac{d_\VV^2}{t}\right).
\end{align*}
For $j\in L$, the integral over $I_j$ is uniformly comparable to 
$(t/a_j)^{\kappa_j+1}$. This gives
\begin{align*}
    D_{t,L}^\kappa(x,y)
    \gtrsim{}&
    t^{-d/2-\sum_{j\in R}\kappa_j}
    \prod_{j\in P}(t+a_j)^{-\kappa_j}
    \prod_{j\in L}a_j^{-\kappa_j}
    \exp\left(-\frac{d_\VV^2}{t}\right).
\end{align*}
Since $a_j\simeq t+a_j$ for $j \in L$ and
$t \simeq t+a_j$ for $j \in R$, this proves the desired lower bound. 

We have proved \eqref{eq:main-simplex-estimate} for $0<t\leq T$. The uniform estimate for $t \geq T$ follows from the
spectral expansion, positivity of the simplex heat semigroup, and the
simple eigenvalue $0$ corresponding to constants, exactly as in the
ball case. This completes the proof of
Theorem~\ref{thm:main-simplex}.

\appendix

\section{Real kernels and holomorphic parameter extensions}
\label{app:holomorphy}

The analytic continuation argument passes through complex parameter values, although the spectral objects in the main theorems are defined for real parameters. We make the distinction precise here. The construction has
five stages:
\begin{enumerate}[label=\textnormal{(\roman*)}]
    \item we start with the reproducing and heat kernels attached to a positive measure in the restricted real parameter range;
    \item we represent these kernels using explicit finite-dimensional formulas that have the correct interpretation also in the extended real parameter range;
    \item we extend these formulas holomorphically in the parameters to a complex domain that, in particular, contains the extended real parameter range;
    \item we find an alternative representation of each kernel in the restricted real parameter range that can also be extended holomorphically to the same complex domain; 
    \item we identify the two extensions by the identity theorem and work in the extended real parameter range using the alternative representation.
\end{enumerate}
Thus, a complex-parameter kernel below is only a holomorphic function
of that parameter. It is not asserted to be the kernel of an
orthogonal projection or a positive heat semigroup.

\subsection{The ball}
We begin with the real object and introduce the standard Jacobi--spherical-harmonic basis on the ball; see \cite[Section~5.2]{DX14}. Let $\mu>-1/2$. For a non-negative integer $n$, the basis of orthogonal polynomials of degree $n$ in $L^2(\dint W_\mu)$ consists of the polynomials
\begin{equation}
\label{eq:ball-explicit-basis}
    Q_{j,\ell}^{n,\mu}(x)
    =
    P_j^{\alpha,\beta_m}(2\norm{x}^2-1)Y_\ell^m(x),
    \qquad
    0\leq j\leq\lfloor n/2\rfloor, \quad 1\leq\ell\leq a_m^d,
\end{equation}
where $m = n-2j$ and $\{Y_\ell^m:1\leq\ell\leq a_m^d\}$ is an orthonormal basis of spherical harmonics of degree $m$ on $\mathbb S^{d-1}$. In \eqref{eq:ball-explicit-basis}, we also used the auxiliary parameters
\[
    \alpha=\mu-\frac12,
    \qquad
    \beta_m=m+\frac{d-2}{2}.
\]
Using polar coordinates
and the Jacobi norm formula, we obtain
\begin{equation}
\label{eq:ball-explicit-norm}
    H_{j,n}^{\mu}
    \coloneqq
    \norm{Q_{j,\ell}^{n,\mu}}_{L^2(W_\mu)}^2
    =
    2^{-\alpha-\beta_m-2}h_j^{\alpha,\beta_m},
\end{equation}
where
\[
    h_j^{\alpha,\beta}
    =
    \frac{2^{\alpha+\beta+1}}{2j+\alpha+\beta+1}
    \frac{\Gamma(j+\alpha+1)\Gamma(j+\beta+1)}
    {j!\Gamma(j+\alpha+\beta+1)}.
\]
Consequently, for real $\mu>-1/2$, we have
\begin{equation}
\label{eq:ball-kernel-finite-sum}
    P_n(W_\mu,x,y)
    =
    \sum_{j=0}^{\lfloor n/2\rfloor}
    \frac{1}{H_{j,n}^{\mu}}
    \sum_{\ell=1}^{a_{n-2j}^d}
    Q_{j,\ell}^{n,\mu}(x)Q_{j,\ell}^{n,\mu}(y).
\end{equation}

We now define the complex extension without introducing a complex orthogonality measure. Note that Jacobi polynomials are polynomials in their
parameters. Indeed, we can write
\begin{equation}
\label{eq:jacobi-complex-polynomial}
    P_j^{\alpha,\beta}(s)
    =
    \frac{1}{j!}\sum_{k=0}^j
    \binom{j}{k}
    (j+\alpha+\beta+1)_k
    (\alpha+k+1)_{j-k}
    \left(\frac{s-1}{2}\right)^k,
\end{equation}
using the Pochhammer symbol; see \cite[Chapter~IV, Section~4.22]{Szego} and \cite[Section~2]{NSS24}. Hence, the right-hand side makes sense for complex $\alpha$ and $\beta$.

We now introduce the complex parameter domain
\[
    \mathcal H_{\BB}
    \coloneqq
    \{\zeta\in\CC:\operatorname{Re}\zeta>-1/2\}.
\]
For $\zeta\in\mathcal H_{\BB}$, let $\alpha_\zeta=\zeta-1/2$ and $\lambda_\zeta = \zeta+\frac{d-1}{2}$. Moreover, define
\[
    \widehat Q_{j,\ell}^{n,\zeta}(x)
    \coloneqq
    P_j^{\alpha_\zeta,\beta_m}(2\norm{x}^2-1)Y_\ell^m(x)
\]
and
\begin{equation}
\label{eq:ball-complex-normalizer}
    \mathcal N_{j,n}(\zeta)
    \coloneqq
    2(2j+\alpha_\zeta+\beta_m+1)
    \frac{j!\Gamma(j+\alpha_\zeta+\beta_m+1)}
    {\Gamma(j+\alpha_\zeta+1)\Gamma(j+\beta_m+1)}.
\end{equation}
For real $\zeta=\mu>-1/2$, we have
$\mathcal N_{j,n}(\mu) = (H_{j,n}^{\mu})^{-1}$. Finally, for arbitrary
$\zeta\in\mathcal H_{\BB}$, we set
\begin{equation}
\label{eq:ball-complex-projector}
    \widehat P_n^\zeta(x,y)
    \coloneqq
    \sum_{j=0}^{\lfloor n/2\rfloor}
    \mathcal N_{j,n}(\zeta)
    \sum_{\ell=1}^{a_{n-2j}^d}
    \widehat Q_{j,\ell}^{n,\zeta}(x)
    \widehat Q_{j,\ell}^{n,\zeta}(y).
\end{equation}
The new object is an explicit finite sum, not
an orthogonal projection for a complex measure.

\begin{proposition}[Holomorphic extension of the ball kernels]
\label{prop:appendix-holomorphy}
Fix $d\geq2$ and $x,y\in\BB^d$.
\begin{enumerate}[label=\textnormal{(\roman*)}]
    \item For every $n\geq0$, the map
    $\zeta\mapsto\widehat P_n^\zeta(x,y)$ is holomorphic in
    $\mathcal H_{\BB}$ and, for real $\mu>-1/2$, 
    \[
        \widehat P_n^\mu(x,y)=P_n(W_\mu,x,y).
    \]
    \item For every $t>0$, the series
    \begin{equation}
    \label{eq:ball-complex-heat-series}
        \widehat h_t^\zeta(x,y)
        \coloneqq
        \sum_{n=0}^{\infty}
        \me^{-tn(n+2\lambda_\zeta)}
        \widehat P_n^\zeta(x,y)
    \end{equation}
    converges locally uniformly in $\mathcal H_{\BB}$ and defines a
    holomorphic function there. Moreover, for real $\mu>-1/2$, it equals the
    genuine ball heat kernel $h_t^\mu(x,y)$.
\end{enumerate}
\end{proposition}

\begin{proof}
Formula \eqref{eq:jacobi-complex-polynomial} shows that every
$\widehat Q_{j,\ell}^{n,\zeta}(x)$ is entire in $\zeta$. If
$\operatorname{Re}\zeta>-1/2$, then all gamma-function arguments in
the numerator of \eqref{eq:ball-complex-normalizer} have positive real
part. The reciprocal gamma factors are entire. Thus
$\mathcal N_{j,n}$, and hence the finite sum
\eqref{eq:ball-complex-projector}, are holomorphic in
$\mathcal H_{\BB}$. For real $\zeta=\mu$, formulas
\eqref{eq:ball-explicit-norm} and
\eqref{eq:ball-kernel-finite-sum} give the asserted equality.

For part~\textnormal{(ii)}, fix a compact set
$K\subset\mathcal H_{\BB}$ and $x,y \in \BB^d$. The formula for Jacobi
polynomials, the addition formula for spherical harmonics, and standard
gamma-ratio estimates give
\[
    \sup_{\zeta\in K}
    \abs{\widehat P_n^\zeta(x,y)}
    \leq
    C_K\exp\big(C_Kn\log(n+2)\big).
\]
At the same time, we have $\abs{\me^{-tn(n+2\lambda_\zeta)}}
    \leq
    \exp(-tn^2+C_Kn)$ for $\zeta \in K$.
Hence, the Weierstrass test proves local uniform convergence of
\eqref{eq:ball-complex-heat-series}, and therefore its sum is
holomorphic. For real $\mu>-1/2$, the summands are exactly those of
\eqref{eq:ball-heat-series}.
\end{proof}

We will write $P_n(W_\zeta,x,y)$ and $h_t^\zeta(x,y)$ for
$\widehat P_n^\zeta(x,y)$ and $\widehat h_t^\zeta(x,y)$. This notation
will cause no confusion: when $\zeta$ is non-real, it
refers only to the holomorphic functions
\eqref{eq:ball-complex-projector} and
\eqref{eq:ball-complex-heat-series}.

\subsection{The simplex}
Let $\Omega_{\VV}=(-1/2,\infty)^{d+1}$. We begin with the real reproducing kernels for the measure $\dint U_\kappa$ with $\kappa\in\Omega_{\VV}$; see \cite[Section~5.3]{DX14}. For $\eta\in\NN_0^d$, put
$\abs\eta=\eta_1+\cdots+\eta_d$ and
\[
    T_j=\sum_{i=j+1}^d\eta_i,
    \qquad
    r_j(x)=1-\sum_{i=1}^{j-1}x_i,
\]
for $1\leq j\leq d$, with empty sums equal to zero. Define
\begin{equation*}
    a_j(\eta,\kappa)
    =2T_j+\sum_{i=j+1}^{d+1}\kappa_i+\frac{d-j-1}{2},
    \qquad
    b_j(\kappa)=\kappa_j-\frac12.
\end{equation*}
By \cite[Proposition~5.3.1]{DX14}, an orthogonal basis
of the degree-$n$ space consists of the polynomials
\begin{equation}
\label{eq:simplex-product-polynomial}
    Q_\eta^\kappa(x)
    =\prod_{j=1}^d
    r_j(x)^{\eta_j}
    P_{\eta_j}^{a_j(\eta,\kappa),b_j(\kappa)}
    \left(\frac{2x_j}{r_j(x)}-1\right)
\end{equation}
with $\abs\eta=n$. This expression is understood as a polynomial also in the case $r_j(x)=0$.
Indeed, applying
\eqref{eq:jacobi-complex-polynomial} reveals that all potentially singular denominator factors cancel.

Next, we set
\[
    S_\kappa=\sum_{i=1}^{d+1}\kappa_i+\frac{d+1}{2}
    =\lambda_\kappa+1,
    \qquad
    c_j(\eta,\kappa)=a_j(\eta,\kappa)+b_j(\kappa)+1,
\]
and define
\begin{equation*}
    w_\kappa^{\VV}
    =\frac{\Gamma(S_\kappa)}
    {\prod_{i=1}^{d+1}\Gamma(\kappa_i+1/2)}.
\end{equation*}
From the Dirichlet integral
\cite[(5.3.3)]{DX14}, it follows that $w_\kappa^{\VV} = (U_\kappa(\VV^d))^{-1}$. The squared norm of
$Q_\eta^\kappa$ with $|\eta|=n$ in the normalized space 
$L^2(w_\kappa^{\VV} \, \dint U_\kappa)$ has reciprocal
\begin{equation}
\label{eq:simplex-inverse-norm}
    \mathcal C_\eta(\kappa)
    =(S_\kappa)_{2n}
    \prod_{j=1}^d
    \frac{\eta_j!}
    {(a_j(\eta,\kappa)+1)_{\eta_j}
     (b_j(\kappa)+1)_{\eta_j}
     (c_j(\eta,\kappa)+\eta_j)_{\eta_j}}.
\end{equation}
This is the norm formula of \cite[Proposition~5.3.1]{DX14} after
cancelling by the polynomial identity
\[
    (c)_{2m}=(c)_m(c+m)_m.
\]
After cancellation, the resulting expression makes sense when $(c)_m=0$.
It agrees with the actual real norm at these parameters by continuity:
on a small compact real parameter neighborhood, the fixed-degree norm integrands are controlled by
$\prod_{i=1}^{d+1}x_i^{\epsilon-1}$ for some $\epsilon>0$, which is
integrable by the Dirichlet formula. The denominator factors in
\eqref{eq:simplex-inverse-norm} are nonzero throughout the real range,
as is also verified below on the larger complex domain.
Consequently, setting $\mathcal E_n = \{ \eta \in \NN_0^d : \abs\eta=n\}$, we obtain that the real projection kernel with respect to $\dint U_\kappa$ is
\begin{equation}
\label{eq:simplex-product-real-projector}
    P_n(U_\kappa,x,y)
    =w_\kappa^{\VV}
    \sum_{\eta \in \mathcal E_n}
    \mathcal C_\eta(\kappa)
    Q_\eta^\kappa(x)Q_\eta^\kappa(y).
\end{equation}
The factor $w_\kappa^{\VV}$ converts the orthonormal-basis kernel for $w_\kappa^{\VV} \, \dint U_\kappa$
to the kernel for $\dint U_\kappa$.

We now introduce the common complex parameter domain
\[
    \mathcal T_{\VV}
    =\{\zeta\in\CC^{d+1}:
      \operatorname{Re}\zeta_i>-1/2,\ 1\leq i\leq d+1\}.
\]
For $\zeta\in\mathcal T_{\VV}$, define
$a_j(\eta,\zeta)$, $b_j(\zeta)$, $c_j(\eta,\zeta)$,
$S_\zeta$, $w_\zeta^{\VV}$, $\mathcal C_\eta(\zeta)$ by the same formulas with $\kappa$ replaced by $\zeta$. Let
$\lambda_\zeta=S_\zeta-1$ and define $Q_\eta^\zeta$ through the
polynomial expansion of \eqref{eq:simplex-product-polynomial}. Finally, set
\begin{equation}
\label{eq:simplex-global-complex-projector}
    \widetilde P_n^\zeta(x,y)
    =w_\zeta^{\VV}
    \sum_{\eta \in \mathcal E_n}
    \mathcal C_\eta(\zeta)
    Q_\eta^\zeta(x)Q_\eta^\zeta(y).
\end{equation}
For non-real $\zeta$, this definition does not require any orthogonal-projection interpretation.

\begin{proposition}[Holomorphic extension of the simplex kernels]
\label{prop:appendix-simplex-holomorphy}
Fix $d\geq2$ and $x,y\in\VV^d$.
\begin{enumerate}[label=\textnormal{(\roman*)}]
    \item For every $n\geq0$, the map
    $\zeta\mapsto\widetilde P_n^\zeta(x,y)$ is holomorphic in
    $\mathcal T_{\VV}$ and, for $\kappa\in\Omega_{\VV}$,
    \[
        \widetilde P_n^\kappa(x,y)=P_n(U_\kappa,x,y).
    \]
    \item For every $t>0$, the series
    \begin{equation}
    \label{eq:simplex-complex-heat-series}
        \widehat H_t^\zeta(x,y)
        \coloneqq
        \sum_{n=0}^{\infty}
        \me^{-tn(n+\lambda_\zeta)}
        \widetilde P_n^\zeta(x,y)
    \end{equation}
    converges locally uniformly in $\mathcal T_{\VV}$ and defines a holomorphic function there. Moreover, for real $\kappa\in\Omega_{\VV}$, it equals the genuine simplex heat kernel $H_t^\kappa(x,y)$.
\end{enumerate}
\end{proposition}

\begin{proof}
Put $\tau_i=\zeta_i+1/2$ and $R_j=\sum_{i=j}^{d+1}\tau_i$. The bases in the denominator of
\eqref{eq:simplex-inverse-norm} satisfy
\[
    a_j(\eta,\zeta)+1=2T_j+R_{j+1},
    \qquad b_j(\zeta)+1=\tau_j,
    \qquad c_j(\eta,\zeta)+\eta_j=2T_j+\eta_j-1+R_j.
\]
On $\mathcal T_{\VV}$, the first two have positive real parts, and so
does the third whenever $\eta_j\geq1$. In the case $\eta_j=0$, the corresponding Pochhammer symbols equal $1$. Thus, every denominator factor is
nonzero on the same set $\mathcal T_{\VV}$ for every degree. Since
$\operatorname{Re}S_\zeta>0$, the gamma factors defining
$w_\zeta^{\VV}$ are also finite and nonzero there.

For completeness, the $j$th factor of $Q_\eta^\zeta$ is the finite
polynomial
\[
    \frac{1}{\eta_j!}\sum_{k=0}^{\eta_j}\binom{\eta_j}{k}
    (c_j(\eta,\zeta)+\eta_j)_k(a_j(\eta,\zeta)+k+1)_{\eta_j-k}
    (x_j-r_j(x))^k r_j(x)^{\eta_j-k}.
\]
It is entire in $\zeta$ and is defined also when $r_j(x)=0$.
Thus, \eqref{eq:simplex-global-complex-projector} is a finite sum of
holomorphic functions on $\mathcal T_{\VV}$. Its real values agree
with \eqref{eq:simplex-product-real-projector}, proving
part~\textnormal{(i)}.

For part~\textnormal{(ii)}, fix a nonempty compact set
$K\subset\mathcal T_{\VV}$ and choose two positive constants
\[
    \delta_K=\min\left(1, \min_{(\zeta, i) \in K \times \{1,\dots, d+1\}}
      \operatorname{Re}\tau_i\right),
    \qquad
    M_K=\max\left(1,\max_{(\zeta, i) \in K \times \{1,\dots, d+1\}}\abs{\tau_i}\right).
\]
Every denominator linear factor in
\eqref{eq:simplex-inverse-norm} has modulus at least $\delta_K$.
Let $L_K=(d+1)M_K+4$. Because $\eta_j+T_j\leq n$, every numerator linear factor in $(S_\zeta)_{2n}$ and in the preceding polynomial expansion has modulus at most $L_K(n+1)$. Moreover, we have
$0\leq x_j\leq r_j(x)\leq1$ on $\VV^d$. Summing the binomial
coefficients in each factor therefore gives
\[
\begin{split}
    \abs{Q_\eta^\zeta(x)}
    \leq (2L_K(n+1))^n (\eta_1! \cdots \eta_d!)^{-1}, \qquad
    \abs{\mathcal C_\eta(\zeta)}
    \leq(L_K(n+1))^{2n} \delta_K^{-3n}
      (\eta_1! \cdots \eta_d!),
\end{split}
\]
uniformly in $(\zeta,x) \in K \times \VV^d$. Set $W_K=\sup_{\zeta\in K}\abs{w_\zeta^{\VV}}$.
There are at most $(n+1)^d$ multiindices of total degree $n$. Hence, for some sufficiently large constant $C_K$, we obtain
\[
\begin{split}
    \sup_{(\zeta,x,y) \in K \times \VV^d \times \VV^d}
    \abs{\widetilde P_n^\zeta(x,y)}
    \leq W_K(n+1)^d
       (4L_K^4\delta_K^{-3}(n+1)^4)^n
    \leq C_K\exp(C_Kn\log(n+2)).
\end{split}
\]
For $n=0$, this simply bounds
$\widetilde P_0^\zeta=w_\zeta^{\VV}$.
Finally, $\lambda_\zeta=S_\zeta-1$ and
$\operatorname{Re}S_\zeta>0$ imply
\[
    \abs{\me^{-tn(n+\lambda_\zeta)}}
    \leq\me^{-tn(n-1)}.
\]
For each fixed $t>0$, the decay in $n$ above dominates the preceding growth bound. Hence, the Weierstrass test proves local uniform convergence of
\eqref{eq:simplex-complex-heat-series}, and therefore its sum is holomorphic. For real parameters $\zeta$, part~\textnormal{(i)} identifies each summand with the corresponding term of \eqref{eq:simplex-heat-series}, proving part~\textnormal{(ii)}.
\end{proof}

As in the ball case, after Proposition~\ref{prop:appendix-simplex-holomorphy}
we write $P_n(U_\zeta,x,y)$ and $H_t^\zeta(x,y)$ also for
$\widetilde P_n^\zeta(x,y)$ and $\widehat H_t^\zeta(x,y)$. For non-real
$\zeta$, these symbols refer to the holomorphic functions above, not to
kernels of a complex-weight Hilbert space.

\subsection{The role of the Dijksma--Koornwinder formula}

The preceding construction supplies the left-hand sides required in the continuation argument independently of any product formula. The
Dijksma--Koornwinder mechanism \cite{DK} supplies the right-hand sides. In the
classical parameter range, the product formula is an ordinary integral
identity against finite positive measures. Following
\cite[Sections~3 and~4]{NSS24}, one subtracts endpoint values so that
the new integrands vanish at every singular endpoint. Proposition
\ref{prop:one-variable-Pi} then shows that the regularized right-hand
side is holomorphic in the larger parameter domain.

For the ball, the classical identity holds for real $\mu\geq0$.
Both its independently constructed left-hand side and its regularized
right-hand side are holomorphic in
$\mathcal H_{\BB}$, so the one-variable identity theorem gives
Propositions~\ref{prop:regularized-reproducing-kernel} and
\ref{prop:regularized-ball-reduction}. For the simplex, the same
argument is applied successively to the coordinates of $\kappa$; the
compatibility established in
Proposition~\ref{prop:appendix-simplex-holomorphy} guarantees that the
left-hand side is the continuation of the original real kernel.

The logical order is therefore as follows. \begin{enumerate} 
\item The classical product formula yields the regularized holomorphic formula. 
\item The real positive-weight kernel yields an independent holomorphic extension. 
\item The identity theorem then ensures that the two extensions coincide. \end{enumerate}
At the end of the argument the parameters are returned to the real
natural range, where the continued functions coincide with the genuine
reproducing and heat kernels.

\section*{Acknowledgments}

The author is deeply grateful to his advisor, Dariusz Kosz, for his
patient verification of the mathematical correctness of this work and
for countless helpful stylistic and editorial suggestions. His careful
reading has considerably improved both the precision of the arguments
and the clarity of their presentation, and his guidance has been
invaluable throughout the author's doctoral studies.

The author would also like to thank Professor Adam Nowak (Institute of
Mathematics, Polish Academy of Sciences) for his contribution to the
author's mathematical development. In particular, the author is
grateful for having his attention drawn to many interesting problems,
including those in the circle of ideas to which this paper belongs,
and for the time Professor Nowak has so generously devoted to him.

\section*{Declaration of generative AI and AI-assisted technologies}

During the preparation of this work the author used generative AI tools, primarily GPT-5.6 Sol
(OpenAI), for
two purposes:
\begin{enumerate}[label=(\roman*)]
    \item to improve the grammar and language of the text;
    \item to check the correctness of some intermediate steps that were needed in order to prove the main results of the paper.
\end{enumerate}
The author takes full responsibility for the final content of the publication.



\end{document}